\documentclass[11pt]{article}
\usepackage{graphicx}
\usepackage{graphicx}  
\usepackage{caption} 
\usepackage{subcaption} 
\usepackage{amscd}
\usepackage{amsmath}
\usepackage{amssymb}
\usepackage{bm}
\usepackage[all]{xy}
\usepackage{float}
\usepackage{comment}

\usepackage[hyphens]{url}
\usepackage{enumitem} 

\usepackage{amsthm}
\numberwithin{equation}{section}
\newtheorem{thm}{Theorem}[section]

\newtheorem{prop}[thm]{Proposition}
\newtheorem{lem}[thm]{Lemma}
\theoremstyle{definition}
\newtheorem{ex}[thm]{Example}
\newtheorem{rem}[thm]{Remark}
\newcommand{\aalp}{\boldsymbol\alpha}
\newcommand{\xxi}{\boldsymbol\xi}
\newcommand{\act}{\lhd \,}
\newcommand{\ad}{\mathrm{ad}}
\newcommand{\alp}{\alpha}
\newcommand{\bbeta}{\boldsymbol\beta}
\newcommand{\pii}{\boldsymbol\pi}

\newcommand{\del}{\delta}
\newcommand{\ddel}{\boldsymbol\delta}
\newcommand{\eee}{\boldsymbol e}
\newcommand{\ep}{\epsilon}
\newcommand{\eeta}{\boldsymbol\eta}

\newcommand{\lam}{\lambda}
\newcommand{\Lam}{\Lambda}

\newcommand{\ome}{\boldsymbol\omega}

\newcommand{\FH}{{\cal F} \!\ltimes \! {\cal H}}
\newcommand{\FHs}{{\cal F} \!\ltimes^{\! s } {\cal H}}
\newcommand{\ggam}{\boldsymbol\gamma}
\newcommand{\gam}{\gamma}

\newcommand{\PPP}{\mathcal P}

\newcommand{\id}{{\mathrm{id}}}

\newcommand{\aaa}{\boldsymbol{a}}

\newcommand{\FFF}{{\cal F}}
\newcommand{\HHH}{{\cal H}}

\newcommand{\QQQ}{{\cal Q}}
\newcommand{\SSS}{{\cal S}}

\newcommand{\R}{{\Bbb R}}

\newcommand{\nab}{\triangledown}

\newcommand{\pr}{\prime}
\newcommand{\prr}{{\prime\prime}}
\newcommand{\sig}{\sigma}

\newcommand{\ssig}{\boldsymbol\sigma}
\newcommand{\rrho}{\boldsymbol{\rho}}

\newcommand{\ttau}{\boldsymbol\tau}

\newcommand{\Ups}{{\cal V}}
\newcommand{\xxx}{\boldsymbol{x}}
\newcommand{\yyy}{\boldsymbol{y}}

\newcommand{\ract}{\leftharpoonup}
\newcommand{\lt}{\! \ltimes \!}
\begin{document}
%
%

%

%
%


\title{Partial Group Symmetry in Figures I:\\
 Semidirect Products and the Six Coins}

\date{}

\author{Takahiro Hayashi \quad}
\maketitle

\begin{abstract}
In this paper, we construct a partial group \(\mathcal{P} (F)\)
that represents the ``partial symmetry'' inherent in a subset \(F\) of \(d\)-dimensional Euclidean space. 
In cases where \(F\) is not connected, \(\mathcal{P} (F)\) captures more detailed information 
than the conventional symmetry group  \(G(F)\). 
To establish a stronger connection between \(\mathcal{P} (F)\) and \(F\), 
we introduce a novel definition of partial group action. 
Furthermore, to characterize \(\mathcal{P} (F)\) in specific cases, we define partial group actions on other partial groups and present a construction of the corresponding semidirect product.
\end{abstract}

\section{Introduction}

Group theory has long served as a fundamental tool 
for describing symmetry in geometric figures and chemical structures. 
However, classical group theory reveals its limitations 
when applied to more complex structures such as Penrose tilings and quasicrystals, 
where the rigid framework of traditional symmetry fails to capture all the nuances of the observed symmetries.

For example, consider the ``Six Coins'' arrangement (see Figure~1). 
While classical group theory might classify this configuration similarly to a rectangle, 
it overlooks the distinct rotational symmetry of each individual coin 
as well as the overall ``\(3\times2\) translation symmetry'' present in the figure 
(refer to Section~5 for our definition of group symmetry).


\begin{figure}[t]
\centering
\includegraphics[width=27mm]{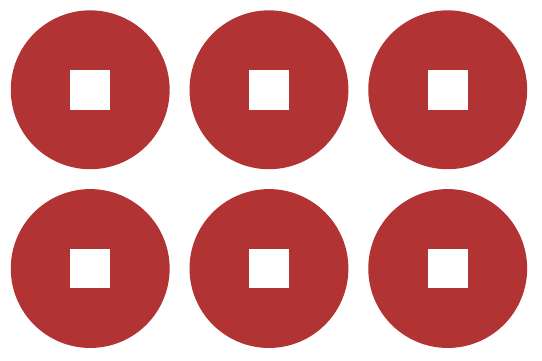}
\caption{Six Coins.}\label{fig1}
\end{figure}

To address these cases of  ``partial symmetry,''
we further develop the framework of symmetry in figures, drawing upon Chermak's introduction of a class of algebraic objects based on partially defined operations  ---  which we refer to as partial groups \cite{Chermak}.
These structures generalize groups (and, more precisely, the associated simplicial sets \cite{Gonzalez}) and were originally identified by Chermak in the context of fusion systems in group theory.
 In fact, Chermak employed partial groups to complete a proof of the Martino-Priddy conjecture, a remarkable problem that lies at the interface of homotopy theory and finite group theory \cite{MP, AKO}.
 Separately, in our own work on constructing generalized Hopf algebras, we arrived at similar concepts of partially defined operations, which further motivated our investigation into their broader properties.
  We believe that the theory of partial groups (and partial monoids) holds promising applications in various fields of 
 mathematics and mathematical science, including 
 quantum measurements  \cite{OCS}.

In this paper, we construct a partial group \(\PPP(F)\) associated with each figure \(F\), 
thereby extending the conventional symmetry group \(G(F)\) to more comprehensively capture the inherent partial symmetry of \(F\). 
More generally, we define a partial group
\[
\PPP(X) = \PPP(X, X_0, G_0)
\]
for any group \(G_0\), a \(G_0\)-set \(X_0\), and a subset \(X \subseteq X_0\). 
 (We note that this construction has already been established in a recent article \cite{OCS}.) 
For instance, when \(G_0\) acts as a transformation group on \(\mathbb{R}^d\) (such as the orthogonal group \(O(d)\) or the Euclidean group \(E(d) \cong \mathbb{R}^d \rtimes O(d)\)), with \(X_0\) representing the set of all figures in \(\mathbb{R}^d\) and \(X\) 
denoting the collection of connected components
of a given figure \(F\), the partial group \(\PPP(F) := \PPP(X)\) 
effectively models the partial symmetry of \(F\).

Our focus lies particularly on figures with two connected components or with mutually parallel, congruent components --- as illustrated by the Six Coins. 
To elucidate the structure of \(\PPP (F)\) for these cases, we develop refined concepts related to partial groups, including sums, fiber sums, and extensions by other groups. 
We also define actions of a partial group on a set and on another partial group
 \(\HHH\). The latter concept allows us to establish a form of semidirect product for partial groups, thereby unveiling new structural insights into \(\PPP(F)\). 
 Our definition of partial group action on \(\HHH\) generalizes that of Grazian and Henke \cite{GH}, yielding a broader formulation of the semidirect product (see also \cite{BG}).

In a forthcoming installment \cite{Hayashi}, we will further explore the partial symmetry properties of semiregular polyhedra.

The structure of this paper is as follows. 
In Section 2, we present a formal definition of partial groups tailored to the objectives of this work, introducing their fundamental structure. Section 3 introduces the concept of partial group action, which plays a crucial role in the subsequent discussions. 
In Section 4, we define the partial group \(\PPP (X)\) and examine its basic properties. 
Section 5 provides a brief overview of the traditional group-theoretical symmetry of figures. 
In Section 6, we initiate an exploration of partial group symmetry in figures, offering an initial analysis. 
Section 7 establishes important connections between specific figures and colimits of groups 
within the category of partial groups.
Section 8 delves into extensions of partial groups by groups, linking these constructions to figures 
with two congruent connected components.
In Section 9, we revisit the simplicial approach to partial groups, 
originally introduced by Broto and Gonzalez \cite{Gonzalez,BG}.
Section 10 details the construction of semidirect product partial groups using the simplicial approach, while Section 11 highlights the universal mapping property of these semidirect products. To this end, we introduce an adjoint (conjugate) action of partial groups on themselves. 
Finally, Section 12 provides detailed descriptions of the partial symmetries of selected figures --- such as the Six Coins --- illustrating the insights gained from the semidirect product construction.

\bigskip
For any set $X$ and  integer $n \geq 1$, we denote its $n$-th Cartesian power 
(i.e., the set of words in $X$ of length $n$)
by $X^n$. 
Given $\xxx = (x_i)_i = (x_1, x_2, \ldots, x_m) \in X^m$ and
$\yyy = (y_j)_j \in X^n$, 
we define their concatenation by 
$\xxx \times \yyy :=(x_1,\ldots, x_m, y_1,  \ldots, y_n) \in X^{m+n}$.
Additionally, we use the following shorthand notation: 
$(x_i)_{p \leq i \leq q}:= (x_p, x_{p+1}, \ldots, x_q)$ 
$(1 \leq p \leq q \leq n)$.
Let $G$ be a group. 
We write $H \leq G$ to indicate that $H$ is a subgroup of $G$. 
For any $c, a \in G$, we define the conjugation by ${}^a c := a c a^{-1}$
and ${c}^{\,a} := a^{-1} c a$.
For a right $G$-set $X$,  the stabilizer of $x \in X$ is denoted by
$G_x$.
We consider the $d$-dimensional Euclidean space
as the set $\R^d$ of the $d$-dimensional real row vectors,
which we treat as a right $E(d)$-set.
For any $x \in \R^d$,  the parallel translation $y \mapsto x + y$  is denoted by 
$t_x \in E(d)$.
Given $F \subseteq \R^d$ and $x \in \R^d$, we define the translated sets by
$F \pm x := \{ y \pm x \mid  y \in F \}$.
For \(n > 0\), we denote the cyclic group of order \(n\) by \(C_n\) 
and the dihedral group of order \(2n\) by \(D_n\).

\section{Partial groups} 
In this section, we introduce the definition of a partial group and related concepts 
that will be used throughout the paper. 
To accurately describe partial group actions, we adopt a notation different from that used in Chermak
\cite{Chermak}.

Let $\PPP_1$ be a non-empty set. For each $n \geq 2$, let $\PPP_n$ be a 
set of sequences $(\pi_1,\pi_2,\ldots, \pi_n)$ of elements of $\PPP_1$, and
let $\nab_n = \nab^\PPP_n\! : \PPP_n \to \PPP_1$ be a map.  
For convenience, we set $\nab_1:= \id_{\PPP_1}$.
We say that $\PPP := ((\PPP_n)_{n \geq 1}, (\nab_n)_{n \geq 1})$
is a {\it  partial semigroup} if the following conditions (P1) and (P2) are satisfied:
\begin{description}[nosep]
\item[(P1)]
If $\pii = (\pi_1,\pi_2,\ldots, \pi_n) \in \PPP_n$ for some $n > 1$, then both
$d_0 ( \pii )= d_0^{\PPP,n} ( \pii ):= (\pi_2,\ldots, \pi_n )$ and
$d_n ( \pii )= d_n^{\PPP,n} ( \pii ):= (\pi_1,\ldots, \pi_{n -1})$ belong to $\PPP_{n-1}$.
\item[(P2)]
For each 
$p, r \geq 0$ and $q > 0$,
\begin{gather}
\label{P2-1}
 \quad (\mathrm{id}_{{\PPP_p}} \times \nab_{q}  \ \times \mathrm{id}_{{\PPP_r}}) (\PPP_{p+q+r}) \subseteq \PPP_{p+1+r},\\
 \label{P2-2}
 \nab_{p+1+r} \circ  (\mathrm{id}_{{\PPP_p}} \times \nab_{q}  \ \times \mathrm{id}_{{\PPP_r}}) |_{ \PPP_{p+q+r}}
=
\nab_{p+q+r}.
\end{gather}
\end{description}

We note that (P1) implies that $\PPP_{p+q+r}$ is a subset of $\PPP_{p} \times \PPP_q \times \PPP_r$.
Hence the left-hand side of \eqref{P2-1} is well-defined.

Let $1 = 1_\PPP$ be an element of $\PPP_1$. 
We say that $1$ is a {\it unit} of $\PPP$, or equivalently, that $\PPP$ is a {\it partial monoid},
if the following two conditions are satisfied:
\begin{description}[nosep]
\item[(P3)]
For each $n > 0$, $\pii = ( \pi_i )_i \in \PPP_n$ and $0 \leq j \leq n$, 
we have $s_j ( \pii ) \in \PPP_{n+1}$, where $s_j = s_j^{\PPP,n}: {\PPP_1}^{\! n }  \to {\PPP_1}^{\! n+1 } $
is defined by 
$s_0 ( \pii ) := (1, \pi_1,\ldots, \pi_n )$, $s_1 ( \pii ) := (\pi_1, 1, \pi_2, \ldots, \pi_n )$, $\dots$,
$s_n ( \pii ) := (\pi_1,\ldots, \pi_n, 1)$.
\item[(P4)]
For each $\pi \in \PPP_1$, we have
\begin{gather}
\label{P4}
\nab_{2} (1, \pi ) = \nab_{ 2} ( \pi, 1 ) = \pi.
\end{gather}
\end{description}

By \cite{Gonzalez}, each partial monoid $\PPP$ becomes a simplicial set (see Section 9).

Having defined a partial monoid, we now extend the structure by introducing an inverse operation. 
We define a \emph{partial group} to be a partial monoid \(\PPP\) equipped with an  operation
$(-)^{-1} : \PPP_1 \to \PPP_1$
(which assigns to each \(\pi \in \PPP_1\) its \emph{inverse} \(\pi^{-1}\)) 
such that the following two conditions are satisfied:
\begin{description}[nosep]
\item[(P5)]
For each $n > 0$ and $\pii = (\pi_1,\pi_2,\ldots, \pi_n) \in \PPP_n$, we have
\begin{gather}
\pii \times \pii^{-1}, 
\pii^{-1} \times \pii \in \PPP_{2n},
\end{gather}
where $\pii^{-1} := (\pi_n^{-1},\ldots,  \pi_2^{-1}, \pi_1^{-1})$.
\item[(P6)]
For each $\pi \in \PPP_1$, we have
\begin{gather}
\nab_{2} (\pi, \pi^{-1})
= 
1
=
\nab_{2} (\pi^{-1},\pi).
\end{gather}
\end{description}

We note that if $(\pi, \rho) \in \PPP_2$ satisfies $\nab_2 (\pi, \rho ) = 1$, 
then 
$\rho= \nab_3 (\pi^{-1}, \pi, \rho ) = \pi^{-1} $ 
by (P5) and (P6).
In particular, we have $(\pi^{-1})^{-1} = \pi$ and 
$\nab_2 (\rho, \sigma)^{-1} = \nab_2 ( \sigma^{-1} , \rho^{-1} )$ 
for each $\pi \in \PPP_1$ and $(\rho, \sigma ) \in \PPP_2$.

Alternatively, we sometimes say that $\PPP_1$ becomes a partial group via
the data $((\PPP_n)_{n \geq 1}, (\nab_n)_{n \geq 1})$, or simply via $\PPP_n$.
When $\PPP_1$ is a finite set, then $\PPP$ is called a {\it finite partial group}.

\par\noindent
\begin{ex}
\label{UpP}
(1)\,\, 
 Every group $G$ can be viewed as a partial group by defining 
$G_n := G^n$ and 
$\nab_n^G (a_1, a_2, \ldots, a_n) := a_1 a_2 \cdots a_n \,\, (a_1, a_2, \ldots, a_n \in G)$.\\
(2)\,\, Let $P$ be a set of points in ${\Bbb R}^d$.
Then 
$\Ups ( {P} )_1: = \{ \overrightarrow{ {\rm P}{\rm Q} } \mid {\rm P}, {\rm Q} \in {P} \} $  
becomes a partial group $\Ups ( {P} )$ via
${\Ups} ( {P} )_n: = 
\{ ( \overrightarrow{ {\rm P}_0 {\rm P}_1 }, 
\overrightarrow{ {\rm P}_1 {\rm P}_2 },  \ldots,
\overrightarrow{ {\rm P}_{n-1} {\rm P}_n })
\mid  {\rm P}_0, {\rm P}_1, \ldots, {\rm P}_n  \in {P} \} $ 
and
$\nab_n 
( \overrightarrow{ {\rm P}_0 {\rm P}_1 }, 
\overrightarrow{ {\rm P}_1 {\rm P}_2 },  \ldots,
\overrightarrow{ {\rm P}_{n-1} {\rm P}_n }) 
:=  \overrightarrow{ {\rm P}_0 {\rm P}_n }.
$
We refer to $\Ups ( {P} )$ as the partial group of {\it vector parades} generated by ${P}$. 
See Section 4 for more general parades." 
For instance, when $d =1$ and $P = \{ 0, 1\}$, then
$\Ups ( {P} ) =\Ups ( \{ 0,1 \} )$ 
coincides with the partial group presented in  \cite[Example 1.2]{ChermakFL1}. 
Specifically, we have
\begin{gather}
\Ups ( \{ 0,1 \} )_n = \{ (\iota_1 - \iota_0, \ldots, \iota_n - \iota_{n-1}) \mid  \iota_0, \ldots, \iota_n \in \{ 0, 1 \}\},\nonumber\\
 \label{V01n}
\nab_n ( \iota_1 - \iota_0, \ldots, \iota_n - \iota_{n-1} ) =  \iota_n - \iota_0.
\end{gather}
In Section 12, we will employ $\Ups ( \{ 0,1 \} )$ to characterize the partial symmetry of figures 
with two congruent or similar connected components.
\end{ex}

Let $\PPP$ and $\QQQ$ be partial groups and 
let
$\phi\!: \QQQ_1 \to \PPP_1$ be a map.
We say that $\phi$  is a 
{\it map of partial groups from $\QQQ$ to $\PPP$}
 if
 $\phi^n (\pii)  \in \PPP_n$
 and
 $\nab_n^{\PPP} (\phi^n ( \pii ))
 =
 \phi ( \nab_n^{\QQQ} ( \pii ))
 $
 for each $n \geq 1$ and $\pii = (\pi_1, \pi_2, \ldots, \pi_n) \in \QQQ_n$,
 where  
 $\phi^n (\pii) = (\phi (\pi_1), \phi (\pi_2), \ldots,\phi (\pi_n))$.
 By abuse of notation, we write $\phi\!: \QQQ \to \PPP$ when
 $\phi$ is a 
 map of partial groups from $\QQQ$ to $\PPP$.
 Similarly to group homomorphisms, 
$\phi$ satisfies  $\phi (1) =1$ and $\phi (\pi)^{-1} = \phi (\pi^{-1})$  
for each $\pi \in \QQQ_1$ (\cite[Section 3]{Chermak}). 
If a bijective map $\phi: \QQQ_1 \to \PPP_1$ gives both a map $\phi : \QQQ \to \PPP$ of partial groups
and a map $\phi^{-1} : \PPP \to \QQQ$ of partial groups, then $\phi$ is called an {\it isomorphism} of partial groups.
For partial groups $\PPP$ and $\QQQ$, we write $\QQQ\cong \PPP$ if there is
an isomorphism $\phi : \QQQ \to \PPP$ of partial groups.
 

Let $\PPP$ be a partial monoid and 
let $\QQQ_1$ be a non-empty subset of $\PPP_1$. For each $n > 1$,
let $\QQQ_n$ be a  subset of ${\QQQ_1}^n \cap \PPP_n$. 
We say that  $\QQQ :=  (\QQQ_n)_{n \geq 1}$ is a 
{\it partial submonoid} 
of $\PPP$ 
(or $\QQQ_1$ becomes a partial submonoid via $  (\QQQ_n)_{n} $)
if the following conditions are satisfied: 
\begin{description}[nosep]
\item[(S1)]
$ \nab_n ( \QQQ_n ) \subseteq \QQQ_1$.
\item[(S2)]
For each $n > 1$, $d_0^{\PPP,n} ( \QQQ_n ), d_n^{\PPP,n} ( \QQQ_n ) \subseteq \QQQ_{n-1}$,
where $d_0^{\PPP,n}$ and  $d_n^{\PPP,n}$ are as in (P1).
\item[(S3)]
For each 
$p, r \geq 0$ and $q > 0$,
$(\mathrm{id}_{{\PPP_p}} \times \nab_{q}  \ \times \mathrm{id}_{{\PPP_r}}) (\QQQ_{p+q+r})
  \subseteq \QQQ_{p+1+r}$.
\item[(S4)]
$1_\PPP \in \QQQ_1$, and 
for each $n > 0$ and
$0 \leq j \leq n$,
$ s_j^{\PPP,n} (\QQQ_{n})
  \subseteq \QQQ_{n+1}$,
where $ s_j^{\PPP,n}$ is as in (P3).  
\end{description}
If, in addition, $\PPP$ is a partial group and $\QQQ$ satisfies:
\begin{description}[nosep]
\item[(S5)]
For each $\pii  \in \QQQ_n$, 
$\pii \times \pii^{-1}, 
\pii^{-1} \times \pii \in \QQQ_{2n}$,
\end{description}
then $\QQQ$ is called a  
{\it partial subgroup} (or an {\it im-partial subgroup}  \cite[Section 9]{ChermakFL3})
of $\PPP$.
A subset $H$ of $\PPP_1$ is called  a {\it subgroup} of $\PPP$
 if it becomes a partial subgroup of $\PPP$ via $(H_n)_n := ({H}^n)_n$
 (\cite{Chermak}). 
It is easy to see that each partial submonoid
(respectively, each partial subgroup) $\QQQ$ of $\PPP$
becomes a partial monoid (respectively, a partial group) 
via ${\nab_n}^{\! \QQQ} := {\nab_n}^{\! \PPP} |_{\QQQ_n}$.
For a map $\phi : \QQQ \to \PPP$ of partial groups, the image $\phi ( \QQQ_1)$
becomes a partial subgroup ${\rm Im} ( \phi )$ of $\PPP$ via
${\rm Im} ( \phi )_n := \{ \phi^n ( \pii ) \mid  \pii \in \QQQ_n \}$.

Although the kernel \( \text{Ker}(\phi) := \{ \pi \in \QQQ_1 \mid \phi(\pi) = 1 \} \) 
also becomes a partial (normal) subgroup \cite{Chermak}, 
it is not necessarily sufficiently large (cf. Remark  \ref{InjKer}). 
Instead of a quotient theory of partial groups by partial normal subgroups \cite{ChermakFL1, Salati}, 
we require a quotient theory based on congruences, 
which is analogous to the quotient theory in semigroup theory 
(see, e.g., \cite{How}).

 Let $\PPP$ be a partial group and
let $\sim$ be an equivalence relation on $\PPP_1$. 
We say that $\sim$ is a {\it (partial group) congruence}  on $\PPP$
if the following two conditions are satisfied: 
\begin{description}[nosep]
\item[(CR1)]
If elements $\pii = ( \pi_i)_i$ and $\rrho = ( \rho_i )_i$ of $\PPP_n$\, $(n>0)$
satisfy
$\pii \sim \rrho$, then $\nab_n ( \pii ) \sim \nab_n ( \rrho )$.
Here,  $\pii \sim \rrho$ 
means $\pi_i \sim  \rho_i$ for each $i$.
\item[(CR2)]
If $\pi,  \rho   \in \PPP_1$
satisfy
$\pi \sim \rho$, then $ \pi^{-1}  \sim \rho^{-1}$.
\end{description}
\begin{prop}
\label{FundT}
{\rm (1)}\,
Let $\sim$ be a congruence on a partial group $\PPP$.
Then $\PPP_1 /\!\!\sim$ becomes a partial group $\PPP /\!\!\sim$ via
\begin{gather*}
 ( \PPP /\!\!\sim )_n 
: = 
\{ ( [ \pi_1 ],  [ \pi_2 ], \ldots,  [ \pi_n ] ) \mid 
( \pi_1, \pi_2, \ldots, \pi_n ) \in \PPP_n \},\\
\nab^{ \PPP /\sim}_n ( [ \pi_1 ],  [ \pi_2 ], \ldots,  [ \pi_n ] )
: =
[ \nab_n^\PPP ( \pi_1, \pi_2, \ldots, \pi_n ) ],
\end{gather*}
where $[ \pi ]:= \{ \rho \mid \rho \sim \pi \}$ for $\pi \in \PPP_1$.\\
{\rm (2)}\,
Let $\phi :\PPP \to \QQQ$ be a map of partial groups.
Then $\phi$  gives rise to a congruence $\sim_\phi$ on $\PPP$ by setting
$ \pi \sim_\phi \rho$ if $\phi ( \pi ) = \phi ( \rho )$.
Moreover, the map $[ \pi ] \mapsto \phi ( \pi )$\,\, $( \pi \in \PPP_1)$  provides an isomorphism 
$\PPP /\!\!\sim_\phi\, \cong\,  {\rm Im} ( \phi )$
of partial groups.
\end{prop}

\begin{proof}
The proof is straightforward.
\end{proof}

\section{Partial group actions}
In this section, we introduce the definition of partial group action on sets, 
which serves as the foundational concept throughout this paper.

Let $\PPP$ be a partial group, and let $X$ be a non-empty set.
For $x \in X$ and $\pii = (\pi_i)_i \in {\PPP_1}^n$, 
we write $(x|\pi_1,\ldots, \pi_n) = (x|\pii )$
 for $(x, \pii) \in X \times {\PPP_1}^n$.
For each $n> 0$, let $X\PPP_{n}$ be a designated subset of $X \times \PPP_n$
and let $\act_{n} = \act_{n}^{X}=\act_{n}^{X \PPP}\!: X\PPP_{n} \to X$,
$(x|\pii) \mapsto \act_n (x|\pii )$
be a map. 
For convenience,  we set
$X \PPP_0 := X$ and
$\act_0 := \mathrm{id}_X$.
We say that 
$X :=  (X, (X\PPP_{n})_{n > 0}, (\act_n)_{n > 0})$ 
is a ({\it right}) $\PPP$-{\it set} --- that is, \(\PPP\) {\it acts} on \(X\)
 (or, equivalently, \(X\) has \(\PPP\)-{\it symmetry}) --- if the following conditions are satisfied:
\begin{description}[nosep]
\item[(A1)] 
If $(x| \pi_1,\ldots, \pi_{n+1}) \in X\PPP_{n+1}$, then 
$(x| \pi_1,\ldots, \pi_{n}) \in X\PPP_{n}$.
\item[(A2)] 
For each $p, q \geq 1$,
\begin{gather}
(\act_{p}\times \id_{\PPP_q} ) (X\PPP_{p+q}) \subseteq X\PPP_{q},\\
\act_{q} \circ (\act_{p} \times \id_{\PPP_q} ) |_{X\PPP_{p+q}}
= \act_{p+q}.
\end{gather}
\item[(A3)] 
For each $p, r \geq 0$ and $q \geq 1$,
\begin{gather}
(\id_{X\PPP_{p}} \times \nab_q \times \id_{\PPP_r}) (X\PPP_{p+q+r}) \subseteq X\PPP_{p+1+r},\\
\act_{p+1+r} \circ ( \id_{X\PPP_{p}} \times \nab_q \times \id_{\PPP_r}) |_{ X\PPP_{p+q+r}} 
= \act_{p + q + r} .
\end{gather}
\item[(A4)] 
For each $n \geq 0$, $\ome = (x | ( \pi_i )_i) \in X \PPP_n$ and $0 \leq j \leq n$,  
we have 
$s_j (\ome)$ $:=$ $(x | \pi_1,\ldots, \pi_j, 1, \pi_{j+1},\ldots,  \pi_n )  \in X \PPP_{n+1}$.
\item[(A5)] 
For each $x \in X$, we have
\begin{gather}
\act_1 ( x | 1 ) = x.
\end{gather}
\item[(A6)] 
If  $(x| \pii) \in  X\PPP_{n} $, then	
$(x| \pii \times \pii^{-1} ) \in X\PPP_{2n}$. 
\end{description}

\medskip
Let $X = (X, (X\PPP_{n})_{n > 0}, (\act^{X}_n)_{n > 0})$
and
$Y = (Y, (Y\PPP_{n})_{n > 0}, (\act^{Y}_n)_{n > 0})$
be $\PPP$-sets.
We say that a map $f:X \to Y$ is a {\it map of $\PPP$-sets} if 
\begin{gather}
(f \times \id_{\PPP_n}) (X\PPP_{n}) \subseteq Y\PPP_{n},\\
\act^{Y}_{n} (f(x)| \pii) 
= 
f ( \act^{X}_{n} (x| \pii) )
\end{gather} 
for each $n>0$ and 
$(x| \pii) \in X\PPP_{n}$.

\begin{ex}
For each partial group $\PPP$, $\PPP_1^{{\rm reg}} := \PPP_1$ and $\PPP_1^{{\rm ad}} := \PPP_1$ 
become right $\PPP$-sets via
$\PPP_1^{{\rm reg}}\PPP_n := \PPP_{1+n}$, $\act^{{\rm reg}}_n := \nab_{1+n}$
and
\begin{align}
\PPP_1^{{\rm ad}}\PPP_{n}
& :=
\{ (\pi_0| \pi_1, \ldots, \pi_n) \in \PPP_1 \times {\PPP_1}^{\! n} \mid  \\
& ({\pi_n}^{\! -1},\ldots, {\pi_1}^{\! -1}, \pi_0, \pi_1, \ldots, \pi_n)
\in \PPP_{1+2n} \},\\
\act^\ad_{n} (\pi_0| \pi_1, \ldots, \pi_n)
& :=
\nab_{1+2n} ({\pi_n}^{\! -1},\ldots, {\pi_1}^{\! -1}, \pi_0, \pi_1, \ldots, \pi_n),
\end{align}
respectively (cf. Lemma \ref{PadP}).
\end{ex}

\begin{ex}
\label{PVP}
Let ${P}$ and $\Ups ({P})$ be as in Example \ref{UpP} (2).
Then ${P}$ becomes a  $\Ups ({P})$-set via
$$
{P}\Ups ({P})_{n} := 
\{ ({\rm P}_0| \overrightarrow{ {\rm P}_0 {\rm P}_1 }, 
\overrightarrow{ {\rm P}_1 {\rm P}_2 },  \ldots,
\overrightarrow{ {\rm P}_{n-1} {\rm P}_n })
\mid {\rm P}_0, {\rm P}_1, \ldots, {\rm P}_n  \in {P} \} ,
$$ 
$$
\act_{n} 
({\rm P}_0|  \overrightarrow{ {\rm P}_0 {\rm P}_1 }, 
\overrightarrow{ {\rm P}_1 {\rm P}_2 },  \ldots,
\overrightarrow{ {\rm P}_{n-1} {\rm P}_n } )
:=  
{\rm P}_n .
$$
\end{ex}

For a partial group $\PPP$ and a $\PPP$-set $X$, we set
$\PPP_{\geq 1} := \coprod_{n \geq 1} \PPP_n$
and
$X \PPP_{\geq 1} := \coprod_{n \geq 1} X \PPP_n$.
For $\pii \in \PPP_n$ and $(x|\pii) \in X\PPP_{n}$,
we write 
$\nab (\pii ): = \nab_n (\pii )$ 
and
$\act (x| \pii ): = \act_{n} (x|\pii )$.

Let $\PPP$ be a partial group and let $X$ be a $\PPP$-set.
We define the $X$-{\it friendly part} $\PPP^{X \text{-fr}}$ of
$\PPP$ to be the partial subgroup of $\PPP$ given by
\[
( \PPP^{X \text{-fr}} )_n 
:= \{ \pii \in \PPP_n \mid  ( x| \pii ) \in X \PPP_n\,\, \text{for some}\,\, x \in X \}.
\]
We say that $X$ is ($\PPP$-){\it friendly}
if $\PPP = \PPP^{X \text{-fr}}$.

\begin{ex}
Let $(\mathcal L,\Delta)$ be an objective partial group in the sense of Chermak \cite{Chermak}. 
Then, $\Delta$  naturally becomes a friendly 
$\mathcal L$-set.
\end{ex}
\section{Parades}
In this section, we introduce what we call the parade construction of a partial group. We anticipate this construction will play a fundamental role in the theory of partial groups, similar to how action groupoids are crucial in groupoid theory.

Let $G_0$ be a group. Let $X_0$ be a right $G_0$-set and
let $X$ be a non-empty subset of $X_0$.
We say that a sequence $\pii = (\pi_1,\pi_2,\ldots, \pi_n)$ of elements of $G_0$ is a {\it parade} in $X$
of {\it length} $n$ if there exists some \(x \in X\) such that
\begin{equation}
\label{parade}
x \pi_1, x \pi_1 \pi_2, \ldots, x \pi_1 \pi_2 \cdots \pi_n \in X.
\end{equation} 
We also say that $\pii$ is a parade with a  {\it starting point} $x$, 
or that $\pii$ defines a {\it route}  $x, x \pi_1, \dots, x \pi_1 \cdots \pi_n$, if condition (4.1) holds.
Let $\PPP (X)_n = \PPP (X, X_0,G_0)_n$ be the set of parades in $X$ of length $n$. 
Define
$\nab_n\! : \PPP (X)_n \to  \PPP (X)_1$
by
$\nab_n (\pi_1,\pi_2,\ldots, \pi_n) := \pi_1 \pi_2 \cdots \pi_n$.
We note that the first part of the following proposition has appeared in a recent article \cite[Example 5.8]{OCS}, and that a special case of this result has already appeared in the original article \cite[Example 2.4 (2)]{Chermak} (see also \cite{Hackney} and \cite{Lynd}).

\begin{prop}
For each $G_0$, $X_0$, and $X$, 
the pair $\PPP (X)$ $=$ $\PPP (X,X_0,G_0)$ $:=$ $ ((\PPP (X)_n)_{n \geq 1},(\nab_n)_{n \geq 1})$
is a partial group.
Moreover, $X$ becomes a friendly right $\PPP (X)$-set by setting
\begin{gather}
 X \PPP(X)_{n} 
 :=
 \{ (x| \pi_1,\ldots, \pi_n )\in X \times \PPP (X)_n \mid 
 x \pi_1, \ldots, x \pi_1 \cdots \pi_n \in X \}, \\
\act_{n} (x| \pi_1,\ldots, \pi_n )
:= x \pi_1 \cdots \pi_n.
\end{gather}
\end{prop}

\begin{proof}
Although the proof is straightforward, we provide it for the first assertion since it is crucial for this paper.
Suppose that $\pii = (\pi_1,\ldots, \pi_n)$ is a parade in $X$ with
a starting point  $x$.
Since
$(\pi_2,\ldots, \pi_n)$ is a parade with a starting point
$y:= x \pi_1 \in X$,  condition (P1) is satisfied.
Next, suppose $n = p + q + r$ and 
$\rho := \pi_{p+1} \pi_{p+2} \cdots \pi_{p+q}$.
Since 
\begin{multline*}
x \pi_1, x \pi_1 \pi_2, \ldots, x \pi_1 \cdots \pi_p,
x \pi_1 \cdots \pi_p \rho,\\
x \pi_1 \cdots \pi_p \rho \pi_{p + q + 1},
\ldots, 
x \pi_1 \cdots \pi_p \rho \pi_{p + q + 1} \cdots \pi_n
\end{multline*}
is a subsequence of the left-hand side of \eqref{parade},
$(\pi_1,\ldots, \pi_p, \rho, \pi_{p + q + 1}, \ldots, \pi_n)$
is a parade with a starting point $x$. Hence (P2) follows from
$$
\nab_{p + 1 + r} (\pi_1,\ldots, \pi_p, \rho, \pi_{p + q + 1}, \ldots, \pi_n)
=
\pi_1 \pi_2 \cdots  \pi_n 
=
\nab_n (\pii).
$$
Set $(\rho_1,\ldots, \rho_n) := \pii^{-1}$ and $z:= x \pi_1 \cdots \pi_n$.
Since the sequence $z, z \rho_1, \ldots$, $z \rho_1 \cdots \rho_n$
is the reverse of $x, x \pi_1, \ldots, x \pi_1 \cdots \pi_n$,
both 
$\pii^{-1} \times \pii$ and $\pii \times \pii^{-1}$ are parades in $X$.
This proves (P5). Since (P3), (P4), and (P6) are obvious, this completes the proof. 
\end{proof}

\par\noindent
\begin{ex}
\label{UpP=P}
(1)\,\,
Suppose that \(P\) is a set of points in \(\Bbb R^d\). 
Let \(X := P\) and let \(X_0 := \Bbb R^d\) be equipped with the translation action of the additive group 
\(G_0 := \Bbb R^d\). 
Then, with these choices, we have $\PPP (X) = \Ups ( {P} )$.\\
(2)\,\, Let $G_0 = G$ be a group and let  $X_0$ be the set $G$ equipped with the adjoint (conjugation) $G$-action.  Let $X = N$ be a (not necessarily normal) subgroup of $G$ .
Then we have $\PPP (X)_n = G^n$ $(n>0)$, since $1^a = 1 \in N$ for every $a \in G$.
The partial group action of $G$ on $N$ is given by 
\begin{gather}
N G_{n} = \{ (m| a_1,\ldots, a_n) \in N \times G^n \mid 
m^{a_1}, m^{a_1 a_2}, \ldots,   m^{a_1 a_2 \cdots a_n} \in N \}, \\
\act (m| a_1,\ldots, a_n) =  m^{a_1 a_2 \cdots a_n}.
\end{gather}
\end{ex}
%
%
%
%

For each $x, y \in X$, set $\PPP (X)_{xy} := \{ a \in G_0 \mid   x a = y \}$. Then we have
\begin{equation}
\label{Pxy}
 \PPP (X )_1  = \bigcup_{x, y \in X}  \PPP (X)_{xy} .
\end{equation} 

\begin{prop}
\label{FinCond}
{\rm (1)}\, For each $x \in X$, the stabilizer $(G_0)_x$ is a subgroup of $\PPP (X)$.\\
{\rm (2)}\,  If $X$ is finite and $(G_0)_x$ is finite for every $x \in X$, then $\PPP (X)$ is finite.
\end{prop}

\begin{proof}
Part (1) is straightforward. 
Suppose $\PPP (X)_{xy}$ contains an element $a$. 
Then the map 
$b \mapsto b a^{-1}$ defines a bijection from 
$\PPP (X)_{xy}$ onto $(G_0)_x $.
Hence, $\PPP (X)_{xy}$ is
finite 
if $(G_0)_x$ is  finite. 
This proves Part (2).
\end{proof}

Let $G_0$ be a group and let $X_0$ be a right $G_0$-set.
For each non-empty subset  $X$ of $X_0$,
we set
$$
 G(X) = G(X, X_0, G_0) := \{ a \in G_0 \mid X a = X \}.
$$
Then $G(X)$ is a subgroup of $G_0$ and $X$ becomes a right 
$G(X)$-set.

\begin{prop}
\label{P=G}
{\rm (1)}\, For every triple $(X, X_0, G_0)$, $G(X)$ is a subgroup of $\PPP (X)$.\\
{\rm (2)}\, If $\PPP (X)_1 = G (X)$ as sets, then  $\PPP (X) = G (X)$ as partial groups; that is,
$\PPP (X) _n = G (X)^n$ for every $n > 0$.\\
{\rm (3)}\,\,
The equality $\PPP (X)_1 = G (X)$ holds if and only if the following two conditions are satisfied:\\
\,\,\,\,{\rm (i)}\, For any \(x, y \in X\), if \(y \in xG_0\) then \(y \in xG(X)\). \\
\,\,\,\,{\rm (ii)}\, For every \(x \in X\), we have \((G_0)_x \subseteq G(X)_x\).
\end{prop}

\begin{proof}
Parts (1) and (2) are straightforward. 
Now, assume that conditions (i) and (ii) hold. 
Let \(\pi \in \PPP(X)_1\). Then there exists \(x \in X\) such that \(y := x\pi \in X\). 
By condition (i), there exists an element \(b \in G(X)\) with \(y = xb\). 
Since \(c := \pi b^{-1}\) belongs to \((G_0)_x \subseteq G(X)_x \subseteq G(X)\), 
it follows that \(\pi = cb \in G(X)\).
Conversely, suppose $\PPP (X)_1 = G (X)$. 
If \(y = x\pi\) for some \(\pi \in G_0\) and \(x, y \in X\), 
then \(\pi \in \PPP(X)_1 = G(X)\), which implies condition (i). 
Moreover, if \(\pi\) is an element of \((G_0)_x\) for some \(x \in X\) (so that \(x\pi = x\)), 
then \(\pi \in \PPP(X)_1 = G(X)\), establishing condition (ii).
\end{proof}

\section{Group symmetry}

In this section, we define the symmetry group $G(F)$ for a figure $F \subseteq \mathbb{R}^d$ and explore its properties under transformations. We also examine the relationship between $G(F)$ and its variants.

Let $F$ be a subset of $\mathbb{R}^d$. 
We define the \textit{symmetry group} $G(F)$, the \textit{special symmetry group} $SG(F)$, 
and the \textit{similarity symmetry group} $SimG(F)$ of the ``figure''$F$ as follows:
\begin{gather*}
G(F)  := G(F,\mathbb{R}^d,E(d)), \quad
SG(F) := G(F,\mathbb{R}^d,M(d)), \\
SimG(F)  := G(F,\mathbb{R}^d,Sim(d)),
\end{gather*}
where $M(d)$ and $Sim(d)$ denote the $d$-dimensional motion group and the group of similarity transformations of $\mathbb{R}^d$, respectively.

\begin{lem}
{\rm (1)} For each $F \subseteq \mathbb{R}^d$, $H \leq Sim(d)$, and $a \in Sim(d)$, we have
\[
G(F,\R^d, a H a^{-1}) = a G( Fa, \R^d, H)a^{-1}.
\]
{\rm (2)} For each $F \subseteq \mathbb{R}^d$ and $a \in E(d)$, the map $b \mapsto a ba^{-1}$ gives an isomorphism $G(Fa) \cong G(F)$ of groups.
\end{lem}

\begin{proof}
This is straightforward.
\end{proof}

Suppose that  $F \subseteq \R^d$ is Lebesgue measurable, and 
 the Lebesgue measure $|F|$ is positive and finite.
 Further, suppose that 
 the coordinate functions $x_i:F\to \Bbb{R}$
 is Lebesgue integrable
 on $F$
for  $1 \leq i \leq d$.  
We define the {\it center of gravity} $z(F) \in \R^d$ of $F$ by
$z(F) := \left(\frac{1}{ |F| } \int_F x_i\, dx_1 \cdots dx_d \right)_i$.
Then, we have $z(Fa) = z(F)a$ for each $a \in E(d)$.

\begin{prop}
\label{G=GO}
Let $F$ and $z := z(F)$ be as above. Then,
\begin{align}
\label{G=GO1}
G(F)
 & = 
 G(F, \R^d,  t_{z}^{\, -1} O(d) t_z) \\
 \label{G=GO2}
 & = 
t_{z}^{\, -1}  G(F-z, \R^d,  O(d)) t_z.
\end{align}
%
\end{prop}

\begin{proof}
Let $a$ be an element of $G(F)$. 
Since $E(d) \cong \R^d \rtimes t_z^{\, -1} O(d) t_z$, 
we have $a = t_x t_z^{\, -1} b\, t_z$ for some $x \in \R^d$ and $b \in O(d)$. 
Substituting this into $z = za$ yields $x = 0$,  which proves the first equality.
The second equality follows from Part (1) of the above lemma.
\end{proof}

For $ F \subseteq \R^d$, we define
its {\it diameter} by $d(F):= \sup \{ |x - y| \mid  x, y \in F \}$. 
Since the diameter \( d(F \gamma a t_x ) \) equals \( \gamma d(F) \) for each \( \gamma > 0 \), \( a \in O(d) \), and \( x \in \mathbb{R}^d \), we obtain the following result:

\begin{prop}
\label{SimG=G}
If
$0 < d(F) < \infty$,
then $SimG(F) = G(F)$. 
\end{prop}

\section{Partial group symmetry}
This section introduces the symmetry partial groups for a figure $F$ in $\mathbb{R}^d$ and examines their properties, 
including conditions under which they coincide with traditional symmetry groups. 

Let $F$ be a subset of $\R^d$.
We define the  {\it symmetry partial group} $\PPP (F)$, 
the {\it special symmetry partial group} $S \PPP (F)$  
and  the  {\it similarity symmetry partial group}
$Sim \PPP (F)$
of the ``figure'' $F$ as follows:
\begin{gather*}
\PPP (F)  := \PPP (cc(F), {\mathfrak P} (\R^d ), E(d)),\quad
S \PPP (F)  := \PPP (cc(F), {\mathfrak P} (\R^d ), M(d)),\\
Sim \PPP (F)  := \PPP (cc(F), {\mathfrak P} (\R^d ), Sim(d)),
\end{gather*}
where $cc(F)$ and  ${\mathfrak P} (\R^d ) = 2^{\R^d}$  denote
the collection of the connected components of $F$ and
the power set of $\R^d$, respectively.

\begin{prop}
\label{PiiF Kihon}
{\rm (1) } For each $F \subseteq \R^d$, we have
\begin{equation}
\label{G = Gcc}
G(F) = G(cc(F), {\mathfrak P} (\R^d ), E(d))
\end{equation}
so that in particular, $G(F)$ is a subgroup of $\PPP (F)$.\\
{\rm (2) }
If $F \subseteq \R^d$ is connected, then as partial groups,  $\PPP (F) = G(F)$. \\
{\rm (3) }
For each $F \subseteq \R^d$ and $ a \in E(d)$, 
the map $\pi \mapsto a \pi a^{-1}$ gives an isomorphism
$\PPP (Fa) \cong \PPP (F)$ of partial groups.\\
{\rm (4) }
For each $K \in cc(F)$, $G(K)$ is a subgroup of $\PPP (F)$.\\
{\rm (5) }
If $|cc(F)|$ is finite and $|G(K)|$ is finite for every $K \in cc(F)$,
then $\PPP (F)$ is finite. 
\end{prop}

\begin{proof}
By continuity, an element \( a \in E(d) \) satisfies \( cc(F)a = cc(F) \) 
if and only if \( Fa = F \). 
This establishes \eqref{G = Gcc}. 
The second assertion of Part (1) follows from Proposition \ref{P=G} (1). 
Parts (2) and (3) are straightforward. 
Finally, Parts (4) and  (5) follow from Proposition \ref{FinCond} (1) and (2), respectively. 
\end{proof}

\begin{figure}[h]
  \centering
  \begin{subfigure}[b]{0.45\textwidth}
    \centering
    \includegraphics[width=33mm]{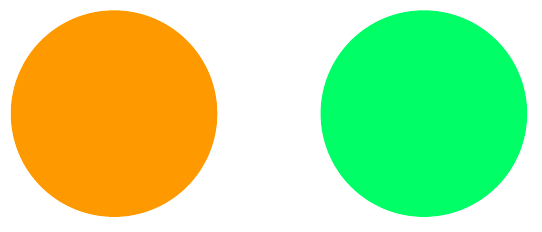}
    \caption{}
    \label{fig:sub-a}
  \end{subfigure}
  \hfill
  \begin{subfigure}[b]{0.45\textwidth}
    \centering
    \includegraphics[width=35mm]{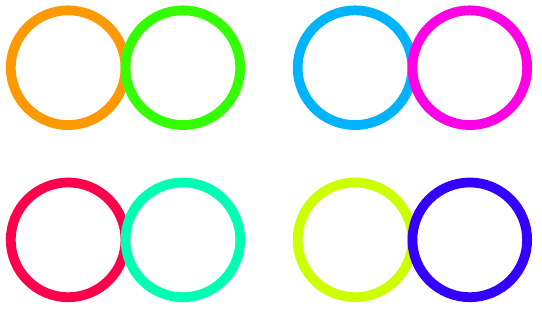}
    \caption{}
    \label{fig:sub-b}
  \end{subfigure}

  \vspace{1em}  
  \begin{subfigure}[b]{0.45\textwidth}
    \centering
    \includegraphics[width=18mm]{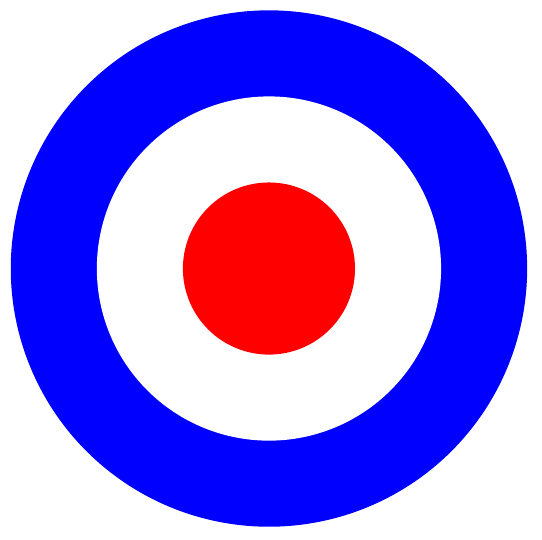}
    \caption{}
    \label{fig:sub-c}
  \end{subfigure}
  \hfill
  \begin{subfigure}[b]{0.45\textwidth}
    \centering
    \includegraphics[width=22mm]{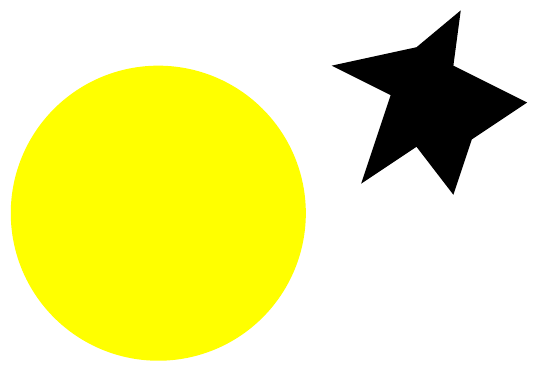}
    \caption{}
    \label{fig:sub-d}
  \end{subfigure}

  \caption{(a) Two disks \(F\). \,\,
    (b) Parametrization of \(P(F)_1\). \,\,
    (c) Disk and annulus \(F'\). \,\,
    (d) Disk with a component without symmetry, \(F''\).}
  \label{fig:overall}
\end{figure}

\begin{ex}
(1)\,\,
Let $F$ be the regular tetragon $[-1,1]^2$ in $\R^2$.
By Proposition \ref{PiiF Kihon} (2) and Proposition \ref{G=GO}, we have  $\PPP (F) \cong G(F, \R^2, O(2)) \cong D_4$.\\
(2)\,\, (Two disks of same size)
Let $F \subset \R^2$ be the union
$D_+ \amalg D_-$ of two disks
$D_\pm  := D \pm e_1$, 
where
$e_1 := (1,0)$ and $D := \{ (x,y) \mid x^2 + y^2 \leq 1/4 \}$. 
(see Figure 2 (a)).
Then, we have:
\begin{gather*}
\PPP (F)_1 
= \left( \PPP (F)_{D_+ D_+} \cup \PPP (F)_{D_- D_-} \right)
\amalg
 \left( \PPP (F)_{D_+ D_-} \cup \PPP (F)_{D_- D_+} \right),\\
%
%
 \PPP (F)_{D_\pm D_\pm,}
= \{ a_{\pm\pm} \mid a \in O(2)\}, 
\quad
 \PPP (F)_{D_\pm D_\mp,}
= \{ a_{\pm\mp} \mid a \in O(2)\}, 
%
\end{gather*}
where $a_{\pm\pm} = t_{\mp e_1} a t_{\pm e_1}$
and $a_{\pm\mp} = t_{\mp e_1} a t_{\mp e_1}$
for each $a\in O(2)$. 
Consequently, the set $\PPP (F)_1$ 
 is parameterized by eight circles, as illustrated in Figure 2 (b). 
 Notice that the four points where these circles intersect correspond precisely to the elements of
$G(F)\cong C_2 \times C_2$.

\par\noindent
(3)\,\,
 Let $F'$ be as depicted in Figure 2 (c). Then, we have $\PPP(F') = G(F') \cong O(2)$,
indicating that the converse of Proposition \ref{PiiF Kihon} (2) does not hold, as
$F'$,
despite not being connected,
satisfies this equality.
The partial group action of $O(2)$ on $cc(F')$ is given by
$cc(F')O(2)_n = cc(F') \times O(2)^n$. 
The action map is trivial, meaning that for every $(K| \pii) \in cc(F')O(2)_n$, we have $\act_{n} (K| \pii) = K$.

\par\noindent
(4)\,\,
Let $F''$ be as depicted in Figure 2 (d).
Then, we have $\PPP(F'') \cong O(2)$ and $G(F'') \cong \{1\}$.
The partial group action of $O(2)$ on $cc(F'')$ is given by
$cc(F'')O(2)_n = \{Y\} \times O(2)^n \cup \{B\} \times \{1\}^n$ with a trivial action map, where $Y$ and $B$ denote the yellow and black components of $F''$, respectively.
Thus, the distinct actions of $O(2)$ on the two-element sets of connected components allow us to distinguish the partial symmetries of $F'$ and $F''$.
\end{ex}

\begin{rem}
\label{PFactF}

In this paper, we primarily regard $\PPP (F)$ and its action on $cc(F)$ as ``invariants'' of the partial symmetry of $F$. 
However, we also expect that $\PPP (F)$ can reveal aspects of the intrinsic structure of $F$, 
since $\PPP (F)$ acts not only on $cc(F)$ but also directly on $F$ itself. This action is given by:
\begin{gather*}
\label{P2-1}
F \PPP (F)_n := \left\{ (x | \pii )\, | \,([x]|\pii ) \in cc(F) \PPP (F)_n \right\},\\
\act_{n} (x| \pi_1,\ldots, \pi_n )
:= x \pi_1 \cdots \pi_n,
\end{gather*}
where $[x]$ denotes the connected component of $F$ containing $x$.
\end{rem}

\section{Wedge sums and fiber sums}

The category of partial groups is cocomplete
(and complete) 
\cite{Chermak, Salati}.
In this section, we provide descriptions of some colimits, which we require.

Let $G_\lam$  be a family of  groups labeled by ${\lam \in \Lambda}$.
For  $a \in G_\lam$  and $b \in G_\mu$,  we write
$a\, {  \sim}\,  b $ if either $\lam = \mu$ and $a = b$, or
$a= 1_{G_{\lam}}$ and $b = 1_{G_{\mu}}$.
Then the quotient  set 
$$
\left({\bigvee}_{\lam \in \Lambda} G_\lam \right)_1
:= 
\left( {\coprod}_{\lam \in \Lambda} G_\lam \right)  /  \sim 
$$
becomes a partial group ${\bigvee}_{\lam \in \Lambda} G_\lam$ via
\begin{gather*}
\left({\bigvee}_{\lam \in \Lambda} G_\lam \right)_n := \bigcup_{\lam \in \Lambda} 
 \{ ([a_1], [a_2],\ldots, [a_n]) \mid 
 a_1, a_2, \ldots, a_n \in G_\lam \},\\
\nab_n ([a_1], [a_2],\ldots, [a_n])
:= [ a_1 a_2 \cdots a_n ]
\quad
(a_1, a_2,\ldots, a_n \in G_\lam),
\end{gather*}
where $[ a ] := \{ b \mid  b \sim a \}$.
Moreover, 
${\bigvee}_{\lam \in \Lambda} G_\lam$ becomes the coproduct of 
$( G_\lam )_{ \lam  \in \Lambda}$
in the category of partial groups through the natural maps
$\iota_\mu : G_\mu \to {\bigvee}_{\lam \in \Lambda} G_\lam$,
$a \mapsto [ a]$.
We call ${\bigvee}_{\lam \in \Lambda} G_\lam$ the {\it wedge sum} of 
$( G_\lam )_{ \lam  \in \Lambda}$.

When each $G_\lam$ is a copy of a group $G$, it is more convenient to define
$\bigvee_\Lam G \cong {\bigvee}_{\lam \in \Lambda} G_\lam$ as follows:
\begin{gather}
\left({\bigvee}_{\!\Lam}\, G \right)_1:= \{ 1_{\bigvee_\Lam G} \}
\coprod
 (\Lam \times (G \setminus \{ 1_G \})),\\
 \left({\bigvee}_{\!\Lam}\, G \right)_n := \bigcup_{\lam \in \Lam} 
 \{ ([\lam, a_1], [\lam,  a_2],\ldots, [\lam, a_n]) \mid 
 a_1, a_2, \ldots, a_n \in G \},\\
\nab_n ([\lam, a_1], [\lam, a_2],\ldots, [\lam,a_n])
:= [\lam, a_1 a_2 \cdots a_n],
\end{gather}
where $[\lam, a] := (\lam, a)$  if $a \ne 1$ and $[\lam, a] := 1_{\bigvee_\Lam G}$ if $a = 1$.

\begin{lem}
Let $G$ be a group and let $(H_\omega)_{ \omega \in  \Omega}$ be a family of groups.
Consider a map of partial groups $\psi: G \to \bigvee_{\omega \in \Omega} H_\omega$, where $\psi(G) \neq \{1 \}$. 
Then there exists a unique 
 $\chi \in \Omega$ and a unique group homomorphism
$f : G \to H_\chi$ such that 
$\psi = \eta_\chi \circ f$,  where $\eta_\chi (h) = [h]$ for $h \in H_\chi$.
\end{lem}

\begin{proof}
Let $a \in G$ be an element such that $\psi (a)  \ne 1$.
By the definition of $\sim$, there exists a unique $\chi \in \Omega$ 
 for which $\psi (a) \in \eta_\chi (H_\chi)$. 
Let $b$ and $c$ be other elements of $G$.
Since $(\psi (a), \psi (b), \psi (c),\psi (bc)) \in ( {\bigvee}_{\omega \in \Omega}\, H_\omega )_4$,
there exist unique $f(b), f(c),f(bc) \in H_\chi$ such that 
$(\psi (b), \psi (c))$ $=$ $(\eta_\chi (f(b)),\eta_\chi (f(c)))$
and that
$\psi (bc)$ $=$ $\eta_\chi (f(bc))$.
Considering $\nab_2 (\psi (b), \psi (c))$ and using the injectivity of $\eta_\chi$,
we obtain $f(b)f(c) = f(bc)$.
\end{proof}

\begin{thm}
Let $(G_\lam)_{ \lam \in  \Lambda}$ and 
$(H_\omega)_{ \omega \in \Omega}$ be families of
nontrivial groups.
Suppose there exists an isomorphism
$$\phi:
{\bigvee}_{\lam \in \Lambda}\, G_\lam 
\cong 
{\bigvee}_{\omega \in \Omega}\, H_\omega
$$
of partial groups.
Then there exist a  bijection
$\varphi\! : \Lambda \cong \Omega$
and group  isomorphisms
$ { f_\lam} \! : G_\lam \cong H_{\varphi (\lam)}$\,\,
$(\lam \in \Lam)$
such that the following diagram is commutative for each $\mu \in \Lambda$:
$$
\begin{CD}
G_\mu @> f_{\mu} >> 
H_{\varphi (\mu )}\\
@ V{\iota_\mu}VV @ VV{\eta_{\varphi (\mu)}}V\\
{\bigvee}_{\lam \in \Lambda}\, G_\lam @>>{\phi} > 
{\bigvee}_{\omega \in \Omega}\, H_\omega
\end{CD}
$$
where $\iota_\mu$ and $\eta_\chi$ are
as above.   
\end{thm}

\begin{proof}
Let $\mu$ be an element of $\Lam$ and let $\chi$ be an element of $\Omega$.
Applying the above lemma to $\psi : = \phi \circ \iota_\mu$, we obtain a unique 
$\varphi (\mu) \in \Omega$ and a unique $f_\mu : G_\mu \to H_{\varphi (\mu)}$ 
such that the above diagram is commutative.
Also, we obtain a unique 
$\varphi^- (\chi) \in \Lam$ and a unique $f_\chi^- : H_\chi \to G_{\varphi^- (\chi)}$ 
such that 
$\iota_{\varphi^- (\chi)} \circ f^-_\chi = \phi^{-1} \circ \eta_\chi$.
Combining these, we get 
$\iota_\mu = \iota_{\varphi^- (\varphi (\mu))} \circ f^-_{\varphi (\mu)} \circ f_\mu$.
Substituting $1 \ne a \in G_\mu$ into this equality, we obtain 
$\varphi^- ( \varphi (\mu))= \mu$. 
This implies the bijectivity of $\varphi$ and $f_\mu$.
\end{proof}

\begin{figure}[h]
\centering
\mbox{\raisebox{0mm}{\includegraphics[width=27mm]{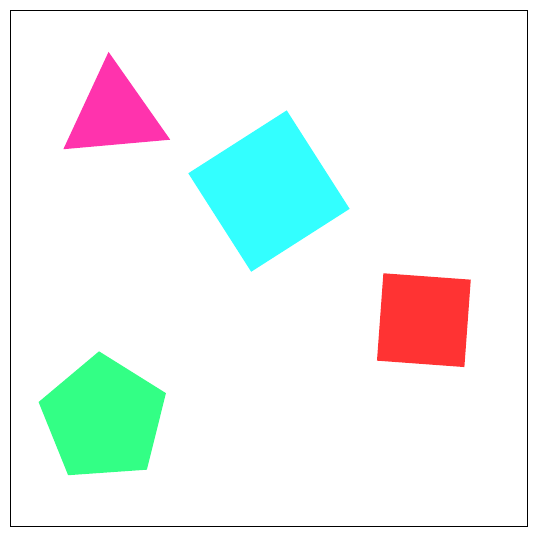}}}
\hspace{20mm}
\mbox{\raisebox{0mm}{\includegraphics[width=27mm]{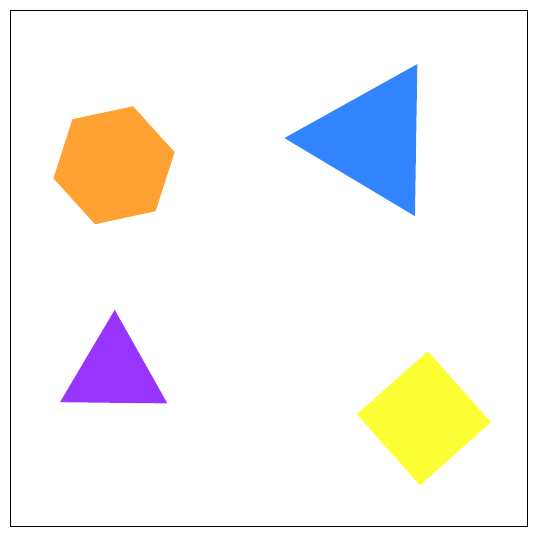}}}
\caption{ Configurations \(F\) and  \(F'\).}
\label{fig3}
\end{figure}

\par\noindent
\begin{ex}
Consider the configurations \(F\) and \(F'\) shown in the left and right panels of Figure~\ref{fig3}, respectively.

\noindent
Then, we observe that $|\PPP (F)_1| = 29 = |\PPP (F')_1|$. 
Therefore, $|\PPP(-)_1|$ cannot distinguish between them. 
However, the theorem above implies that
$\PPP (F) \cong D_3 \vee D_4 \vee D_4 \vee D_5  \not\cong 
D_3 \vee D_3 \vee D_4 \vee D_6 \cong
\PPP (F')$.
\end{ex}

Let $G_+$, $G_-$, and $H$ be groups, 
and let $f_\pm : H \to G_\pm$ be group homomorphisms. 
We denote the pushout (or {\it fiber sum}) of $f_\pm : H \to G_\pm$ 
in the category of partial groups by $ G_+ \vee_H G_-$.

\begin{rem}
\label{InjKer}
Let \( \psi : G_+ \vee G_- \to G_+ \vee_H G_- \) 
be the map induced by the universal mapping property of the wedge sum \( G_+ \vee G_- \). 
Although \( \psi \) is not injective unless \( H = \{1\} \), 
we have \( \psi(\pi) = 1 \) only when \( \pi = 1 \).
 In particular, \( G_+ \vee_H G_- = \text{Im}(\psi) \) 
 is not isomorphic to \( (G_+ \vee G_-) / \mathrm{Ker}(\psi) \) if \( H \neq \{1\} \), 
 where the latter is described in  \cite{Salati}.
Let \( \sim_\psi \) be as in Proposition 2.2 (2). 
Then, for \( \pi, \rho \in (G_+ \vee G_-)_1 \), 
we have \( \pi \sim_\psi \rho \) 
if and only if either \( \pi = \rho \)
or \( \pi = \iota_\pm(f_\pm(h)) \)
and \( \rho = \iota_\mp(f_\mp(h)) \) for some  \( h \in H \).
\end{rem}

Next, suppose that $G_+$ and $G_-$ are subgroups of a group $G_0$.
It is straightforward to verify that $G_+ \cup G_-$ becomes a partial subgroup
of $G_0$ via  
$(G_+ \cup G_-)_n := {G_+}^n \cup {G_-}^n$.
By abuse of notation, we denote this partial subgroup simply by
$G_+ \cup G_-$. 

\begin{prop}
\label{GCupG}
As partial groups, we have $G_+ \cup G_- \cong G_+ \vee_H G_-$, where
$H := G_+ \cap G_-$.
\end{prop}

\begin{proof}
Let $\phi_\pm: G_\pm \to \PPP$ be maps of partial groups that satisfy
$\phi_+ (h) = \phi_- (h)$ for each $h \in H$. 
By gluing \(\phi_+\) and \(\phi_-\) together, we obtain a map
\((G_+ \cup G_-)_1 \to \PPP_1\), which extends to a map
\(\phi : G_+ \cup G_- \to \PPP\) of partial groups.
It is easy to see that $G_+ \cup G_-$ is the pushout of 
$H \hookrightarrow G_\pm$ via the assignment $(\phi_+, \phi_-) \mapsto \phi$. 
\end{proof}

\begin{prop}
Let $F$ be a subset of $\R^d$ consisting of 
two connected components $F_\pm$. 
Then, the inclusions $G(F_\pm) \to \PPP (F)$ induce an injective map 
$G(F_+) \vee_{G(F_+) \cap G(F_-)} G(F_-) \to \PPP (F)$ of partial groups. 
This map is an isomorphism if and only if
$F_+$ and $F_-$ are not congruent.
\end{prop}

\begin{proof}
This is an immediate consequence of the above proposition.
\end{proof}

\begin{figure}[h]
\centering
    \includegraphics[width=0.21\linewidth]{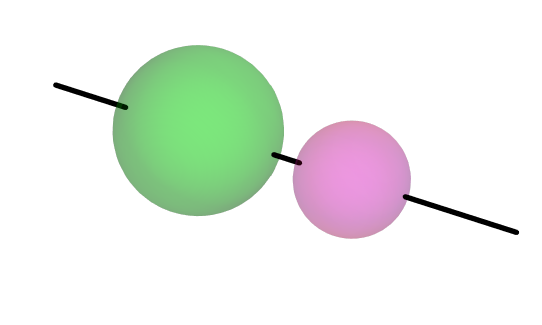}
\caption{Two balls of different sizes.}\label{fig4}
\end{figure}

\begin{ex}
Let \(F\) be the disjoint union of two \(d\)-dimensional balls of different sizes, namely, 
\(B_{r_+} + e_d\) and \(B_{r_-} - e_d\), where 
$B_{r_\pm} := \{ x \in \mathbb{R}^d \mid |x| \le r_\pm \}$,
\(0 < r_- < r_+ < 1\), and
$e_d := (0,\ldots,0,1) \in \mathbb{R}^d$ (cf. Figure 4).
Then, we have  
\begin{align*} 
\PPP (F) & = O(d)_+ \cup O(d)_-\\
           & \cong O(d) \vee_{O(d-1)} O(d),
\end{align*}
where $O(d)_\pm := t_{\pm e_d}^{\, -1} O(d) t_{\pm  e_d}$.
\end{ex}

\section{Extension of partial groups by groups}
In this section, we consider extensions of partial groups $\PPP$ by a group $G$, 
and explore their relation to certain geometric figures. 
This extension theory is a specific case of the framework developed by Broto and Gonzalez  \cite{BG}.

Let $G$ be a group, and
let $\PPP$ be a partial group 
provided with a right $G$-action
$(\pi, g) \to \pi^g$\, $(\pi \in \PPP_1, g \in G)$ on $\PPP_1$.
We say that $ \PPP $ is a ({\it right}) {\it $G$-partial group} if, for every $g \in G$,  
the map 
$\pi \mapsto \pi^g \,\, (\pi \in \PPP_1)$
is an automorphism of $\PPP$ (\cite[Remark 3.2]{GH}).

\begin{ex}
\label{GActVee}
\rm{(1)}\,
Let $G$ and $H$ be groups.
Let $X$ be a right $G$-set. Then 
${\bigvee}_{X}\, H$ becomes a $G$-partial group via $[x,a]^g := [xg, a]$.\\
\rm{(2)}\, Let $H_\lam$ be a family of groups indexed by  $\lam \in \Lam$.
Suppose each $H_\lam$ has a  $G$-action that is compatible with its group structure. 
Then
${\bigvee}_{\lam \in \Lam}\, H_\lam$ becomes a $G$-partial group via 
$[a]^g := [a^g]$\,\, $(a \in H_\lam, g \in G)$.
\end{ex}
For a partial group $\PPP$, we denote its {\it normalizer} by $N( \PPP )$  (\cite{Gonzalez, BG}).
By \cite[Lemma 2.9 (a)]{BG}, $N( \PPP )$ is 
the set of elements  $g$ of $\PPP_1$
satisfying
%
\begin{gather*}
(\pi_1, , \ldots, \pi_{i}, g^{\pm 1}, g^{\mp 1}, \pi_{i+1}, \ldots, \pi_n) \in \PPP_{n+2} 
\end{gather*}
%
for each  $1 \leq i < n$ and $(\pi_1, \ldots, \pi_{n}) \in \PPP_n$.
\begin{prop} 
\label{NP} 
{\rm  (1)}\,The set $N( \PPP )$ is a subgroup of $\PPP$.\\
{\rm (2)}\,  Moreover, $\PPP_1$ becomes a two-sided  $N( \PPP )$-set  via
 \begin{gather}
\label{nualpnu}
 g\cdot \pi \cdot h:= \nab_3 (g, \pi, h) 
 \quad (\pi \in \PPP_1, g, h \in N (\PPP ) ).
\end{gather}
 {\rm (3)}\,  Also, $\PPP$ becomes a right $N( \PPP )$-partial group  via
\begin{gather}
\label{alpnu}
 \pi^g:= g^{-1}\cdot  \pi \cdot g 
 \quad (\pi \in \PPP_1, g \in N (\PPP ) ).
\end{gather}
 {\rm (4)}\, For each $G_0$-set $X_0$ and $X \subseteq X_0$,
 $G(X)$ is a subgroup of $N(\PPP (X))$.\\
 {\rm (5)}\, For each group $G_0$ and subgroups $G_\pm \leq G_0$,
 the intersection $G_+ \cap G_-$ is a subgroup of $N(G_+ \cup G_-)$.
\end{prop}

\begin{proof}
Parts (1) and (3) are established  in \cite{BG}. 
Let $g$ be an element of $N(\PPP)$, and let $(\pi_j)_j$ be an element of $\PPP_n$. 
Since $(1, g^{-1}, g, \pi_1, \ldots, \pi_n) \in \PPP_{\geq 1}$,  it follows that
$(g, \pi_1, \ldots, \pi_n) \in \PPP_{n+1}$ and also
$(\pi_1, \ldots, \pi_n, g) \in \PPP_{n+1}$.
Hence, for each $g, h \in N ( \PPP )$, we have 
$(g, h,\pi)$, $(g, \pi, h), (\pi, g, h)\in \PPP_3$.
Part (2) follows immediately from these observations.
Parts (4) and (5) are straightforward.
\end{proof}

%

Let $\HHH$ be a partial group, and let $G$ be a group.
Let 
$\ract : \HHH_1 \times G \to \HHH_1$, $(\eta, a) \mapsto \eta \ract a$ 
and $\sig : G \times G \to N (\HHH )$ be maps.
We say that \((\ract, \sig)\) is a \textit{factor set} for \((\HHH, G)\)
 (or a \(G\)-\textit{twisting pair} for \(\HHH\); 
see \cite[Definition 4.4]{BG}) if the following conditions are satisfied:
\begin{description}[nosep]
\item[(FS1)]
For each $\eta \in \HHH_1$ and $a \in G$, $\eta \ract 1 = \eta$ and
    $\sig (1, a) = \sig (a,1) = 1$.
\item[(FS2)] 
For each $\eeta = (\eta_i)_i \in \HHH_{\geq 1}$ and $a \in G$, 
we have
$\eeta \ract a \in \HHH_{\geq 1}$
and $\nab (\eeta \ract a ) = \nab (\eeta) \ract a$, 
where 
$\eeta \ract a := (\eta_i \ract a)_i$.
    \item[(FS3)] For each $\eta \in \HHH_{1}$ and $a, b \in G$, 
$(\eta \ract a) \ract b$ $=$ $(\eta \ract ab)^{\sig (a,b)}$.
    \item[(FS4)] For each  $a, b, c \in G$, 
$\sig (a,bc) \sig (b,c)$ $=$ $\sig (ab,c) (\sig (a,b)\ract c)$.
\end{description}

\begin{lem}
Let $(\ract, \sig )$ be a {\it factor set } 
for $( \HHH, G)$. \\
\rm{(1)}\,
For each 
$a \in G$, the map $\eta  \mapsto \eta \ract a$
defines an
automorphism  of the partial group $\HHH$.
In particular, 
$(\eta \ract a)^{-1} = \eta^{-1} \ract a$ for each $\eta \in \HHH_1$.\\
\rm{(2)}\,
For each $a \in G$ and $g \in N( \HHH )$, we have $  g \ract_N a:=  g \ract a \in N(\HHH)$. 
Moreover, $(\ract_N, \sig )$ forms a factor set for $(N(\HHH), G)$.  
\end{lem}

\begin{proof}
Part (1) follows from
Proposition \ref{NP} (3) by setting $b=a^{-1}$ in (FS3).
Let $(\eta_j )_j$ be an element of $\HHH_n$.
Since 
\[(\eta_1 \ract a^{-1}, \ldots, \eta_i \ract a^{-1}, g^{\pm 1}, g^{\mp 1}, 
\eta_{i+1} \ract a^{-1}, \ldots, \eta_n \ract a^{-1}) \in \HHH_{n + 2}, 
\]
we have 
$({\eta_j}^{\sig (a^{-1},a)})_{1 \leq j \leq i} \times
((g \ract a)^{\pm 1}, (g  \ract a)^{\mp 1}) \times 
({\eta_j}^{\sig (a^{-1},a)})_{i < j \leq n} \in \HHH_{n + 2}$.
This confirms that $g  \ract a \in N(\HHH)$. 
The remaining assertion in Part (2) is straightforward.
\end{proof}

Let $(\ract, \sig )$ be a factor set for $( \HHH, G)$.
For each $n > 0$, we set
 \[
( G \ltimes_{\ract,\sig } \HHH )_n := 
\{ \aalp \in ( G \times \HHH_1 )^n \mid 
A( \aalp ) \in \HHH_n \},
\]
where
\begin{gather}
A ( \aalp )
:=
( \eta_1 \ract (c_2, c_3, \ldots c_n),
\ldots,
\eta_{n-1} \ract c_n  , \eta_n ),\\
\eta_i \ract (c_{i+1}, c_{i+2}, \ldots c_n) := (\cdots (\eta_i \ract c_{i+1})\ract  \cdots 
\ract c_{n-1}) \ract c_n
\end{gather}
for each $\aalp = ((c_i, \eta_i))_i \in (G \times \HHH_1)^n$.
By Broto and Gonzalez \cite[Theorem 4.5]{BG}, 
$G \ltimes_{\ract,\sig } \HHH$  $:=$
$ (( G \ltimes_{\ract,\sig } \HHH )_n)_{n \geq 1}, (\nab_n^{ G \ltimes_{\ract,\sig } \HHH })_{n \geq 1})$
becomes a partial group, where
\begin{gather}
\nab^{G \ltimes_{\ract,\sig } \HHH} ( \aalp )
:= ( c_1 c_2 \cdots c_n,  \nab^{\HHH}_{n+1} (\sig (c_1,\ldots,c_n), A ( \aalp ))),\\
\sig (c_1,c_2,\ldots, c_n):= \sig (c_1, c_2 \cdots c_n)\,\sig (c_2, c_3 \cdots c_n) \cdots \sig (c_{n-1},c_n).
\end{gather}
The unit of $G \ltimes_{\ract,\sig } \HHH$ is
$( 1_G, 1_\HHH )$.
When $\sig \equiv 1$, $\HHH$ becomes a $G$-partial group 
via $\eta^c := \eta \ract c$.
In this case, we write $G \ltimes \HHH$ for $G \ltimes_{\ract,\sig } \HHH$.

\begin{ex}
Consider two configurations of the set \(F = F_+\amalg F_-\): the first is shown in the left panel of Figure~\ref{fig4}, and the second in the right panel of Figure~\ref{fig4}.  
In addition to the \(C_2\)-symmetry given by reflection across the black line, \(F\) possesses \(C_3\vee C_3\)-symmetry.  
Combining these, one sees that \(F\) admits a \(C_2\ltimes (C_3\vee C_3)\)-action.  
The \(C_2\)-action on \(C_3\vee C_3\) is described by Example~\ref{GActVee}(1) for the configuration depicted in the left panel, and by Example~\ref{GActVee}(2) for the configuration depicted in the right panel.
\end{ex}

\begin{figure}[t]
\centering
    \includegraphics[width=30mm]{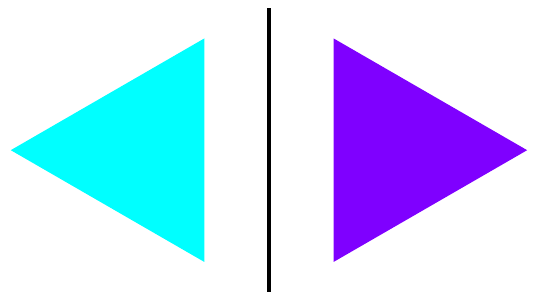}
\hspace{20mm}
    \includegraphics[width=36mm]{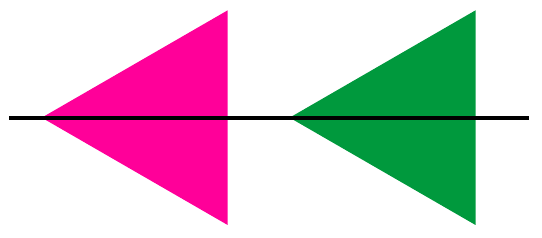}
\caption{Reflection symmetry of two triangles.}
\label{fig4}
\end{figure}

\begin{thm}
\label{P=C2GG} 
Let $F$ be a subset of $\R^d$  with two connected components, $F_+$ and $F_-$.
Assume there exists an isometry $a \in E( d )$ such that $F_+ a = F_-$  
and $F_- a = F_+$.
Then, 
\[
 \PPP  (F) \cong 
C_2 \ltimes_{\ract,\sig } \left( G(F_+) \vee_{G(F_+) \cap G(F_-)} G(F_-)\right)
\]
for some $(\ract, \sig )$.
Identifying $G(F_+) \vee_{G(F_+) \cap G(F_-)} G(F_-)$
with $\HHH := G(F_+) \cup G(F_-)$ via Proposition \ref{GCupG},
the factor set $(\ract, \sig )$ is defined by
$\eta \ract c := a^{-1} \eta a$ and $\sig (c,c) := a^2$,
where $c$ is the generator of $C_2$. 
In particular, $\PPP (F) \cong C_2 \ltimes \HHH$ if $a^2 = 1$.
If $a \in M(d)$, we also have
\[
 S\PPP  (F) \cong 
C_2 \ltimes_{\ract,\sig } \left( SG(F_+) \vee_{SG(F_+) \cap SG(F_-)} SG(F_-)\right).
\]
Explicitly, the inverse $\psi$ of the first isomorphism  is given by
$\psi (c^\iota, \eta) = a^\iota \eta$
for \(\iota \in \{0,1\}\) and \(\eta \in \HHH_1\).
\end{thm}
\begin{proof}
By Proposition \ref{NP} (5), we have $a^2 \in N(\HHH)$.
It is straightforward to verify that $(\ract, \sig )$ is a factor set.
To show that the above map $\psi$ is indeed a map of partial groups,
let $\aalp =((c^{\delta_i},\eta_i))_i$ be an element of  
$(C_2 \ltimes_{\ract,\sig } \HHH)_n$,
where $\del_1, \ldots, \del_n \in \{ 0, 1 \}$. 
Define $(\rho_i)_i := A (\aalp ) \in \HHH_n$.
Then, we have 
\begin{equation}
\label{rhoi}
 a^{- \del_{i+1} - \cdots -\del_n} \eta_i a^{\del_{i+1} + \cdots +\del_n} = \rho_i  \in G(F_\ep )
\quad (1 \leq i \leq n )
\end{equation}
for some $\ep \in \{ \pm 1\}$, where $F_{\pm 1} := F_\pm$. 
Hence,
\[
F_{\ep_{i-1}} a^{\del_i}  \eta_i 
= 
F_\ep \rho_i a^{-\del_{i+1} - \cdots -\del_n} 
= F_{\ep_i},
\]
where $\ep_i := \ep (-1)^{\del_{i+1} + \cdots +\del_n}$
for $0 \leq i \leq n$,
This implies $\psi^n (\aalp )  \in \PPP (F)_n$.
Using induction on $n$, we obtain
$\sig (c^{\del_1}, c^{\del_2}, \ldots, c^{\del_n}) = a^{\del - \iota}$,
where 
$\iota = 0$ if $\del := \del_{1} + \cdots +\del_n$ is even and $\iota = 1$ if $\del$ is odd.
Applying this result with \eqref{rhoi}, we find
$\psi ( \nab ( \aalp )) = a^\iota a^{\del - \iota} \rho_1 \cdots \rho_n = \nab ( \psi^n ( \aalp ))$. 
Next, suppose $\pii = (\pi_i)_i  \in \PPP (F)_n$.
Then, we have $F_{\ep_{i-1}} \pi_i  = F_{\ep_i}$
for some $\ep_i \in \{ \pm 1\}$\,\, $(0 \leq i \leq n)$.
This implies $\eta_i := a^{- \del_i} \pi_i \in G(F_{\ep_i})$, 
and also $A (((c^{\delta_i},\eta_i))_i) \in G( F_{\ep_n })^n$,
where 
$\del_i = 0$ if $\ep_{i-1} = \ep_i$ and  
$\del_i = 1$ if $\ep_{i-1} \ne \ep_i$. 
This proves the bijectivity of $(C_2 \ltimes_{\ract,\sig } \HHH)_n \to \PPP (F)_n$, $\aalp \mapsto \psi^n (\aalp )$.
\end{proof}

\begin{figure}[h]
\centering
    \includegraphics[width=0.21\linewidth]{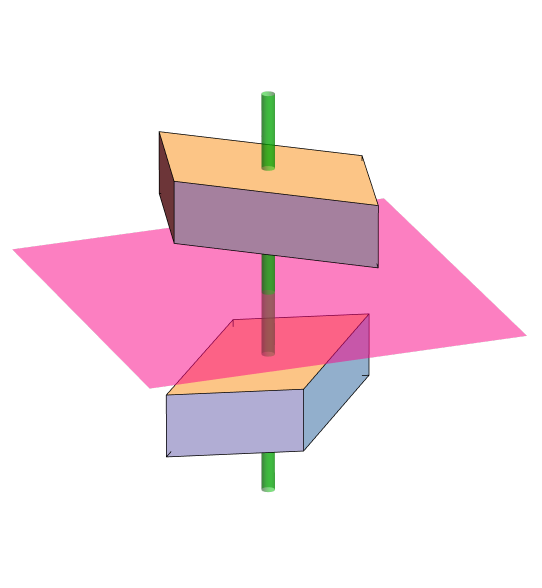}
\caption{Two prisms.}\label{fig6}
\end{figure}

\begin{ex}
Let \(a \in O(3)\) be the composition of a \(90^\circ\) rotation about the \(z\)-axis 
with a reflection across the \(xy\)-plane. 
Consider prisms  $F_+$ and $F_- := F_+ a$ with empty intersections such that 
$z ( F_\pm ) = (0,0,\pm 1)$ and that 
their bases are parallelograms that are parallel to
the $xy$-plane (cf. Figure 6).
Since \(\sig (c,c) = a^2 = \mathrm{diag}(-1,-1,1) \neq 1\), the figure 
$F := F_+ \coprod F_-$
provides an example of the above theorem with a non-trivial cocycle.
Unless the bases of $F_\pm$ are 
 rhombuses, there exists no isometry $\tilde{a}$ such that 
 $F_\pm \tilde{a} = F_\mp$ and that ${\tilde{a}}^2 = 1$.
\end{ex}

\section{Prepartial group}
In this section, we revisit the simplicial definition of a partial group as introduced by Broto and Gonzalez  \cite{BG}. 
For clarity, we distinguish their version from the original definition, 
referring to it as the  ``prepartial group.''
The results of this section will be used in the next section to construct semidirect products of partial groups.

Let \( \Delta \) be the full subcategory of the category of ordered sets, where the objects are 
\(\text{ob}\Delta = \{[n] \mid n \geq 0 \}\) with \([n] := \{0,1,\dots,n\}\).
A {\it simplicial object} in a category \( \bf C \) is a contravariant functor from \( \Delta \) to \( \bf C \).
Specifically, a {\it simplicial set} is a simplicial object in the category of sets.
We say that a simplicial set $\SSS$ is {\it reduced} if  \( \SSS_0 \) contains a unique element \( * \),
where \(\SSS_n := \SSS([n])\) for $n \geq 0$. 
For a simplicial set \(\SSS\), 
the {\it face operators} \(d_i = d_i^{n} = d_i^{\SSS,n}: \SSS_n \to \SSS_{n-1}\) 
and  {\it degeneracy operators} 
\(s_i  = s_i^{n} = s_i^{\SSS,n}: \SSS_n \to \SSS_{n+1}\) (for \(0 \leq i \leq n\)) 
are defined by \(d_i := \SSS(d^i)\) and \(s_i := \SSS(s^i)\). 
The maps \(d^i\) and \(s^i\) are given by \(d^i(j) = j\) if \(j < i\), and \(d^i(j) = j + 1\) 
if \(j \geq i\); similarly, \(s^i(j) = j\) if \(j \leq i\), and \(s^i(j) = j - 1\) if \(j > i\).
 Following Broto and Gonzalez \cite{BG}, we define 
the {\it spine operator} 
(also known as the {\it (Grothendieck-)Segal map})
$\eee^n : \SSS_n \to (\SSS_1)^n$ for $n \geq 1$ by 
$\eee^n (\ssig) = (e_i^n (\ssig))_i$, where
 $e_i^n = d_0^{\, i, i - 1} \circ d_n^{\,n, n - i}$, with 
 $d_0^{\, n, i } := d_0^{\, n - i + 1} \circ \cdots \circ d_{0}^{\, n-1} \circ d_0^{\, n}$,
and
 $d_n^{\, n, i } := d_{\,  n-i + 1}^{\, n-i + 1} \circ \cdots \circ d_{n - 1}^{\, n-1} \circ d_n^{\, n}$.

We define a {\it prepartial monoid}  (also referred to as an
{\it edgy reduced simplicial set} \cite{HL})
\( \SSS \) 
as a reduced simplicial set 
where the spine operator \( \eee^n \) is injective for each \( n \). 
An  {\it  inversion} for a prepartial monoid \( \SSS \) is 
a sequence \( \nu = (\nu_n)_{n \geq 0} \)
of maps \(\nu_n : \SSS_n \to \SSS_n\) that satisfy the following conditions:

\begin{description}[nosep]
\item[(Inv1)]
     \( \nu_n \circ \nu_n = \text{id} \), $d_i^n \circ \nu_n = \nu_{n-1} \circ d_{n - i}^{n}$, 
    and $s_i^n \circ \nu_n = \nu_{n+1} \circ s_{n - i}^{n}$  for \( 1 \leq i \leq n \).
\item[(Inv2)]
 For each \( n \geq 1 \) and $\ssig \in \SSS_n$, 
    there exists  $\ttau \in \SSS_{2n}$ such that
    $d_{2n}^{\, 2n,n} ( \ttau) = \nu_n ( \ssig )$,
    $d_{0}^{\, 2n,n} ( \ttau) = \ssig$,
    and 
    $d_1^{\, 2}  \circ \cdots \circ d_1^{\,2n} (\ttau ) = s_0 (*)$.
\end{description}

A  {\it prepartial group} is a prepartial monoid that is equipped with such an inversion.
According to \cite{Gonzalez},  each partial monoid \( \PPP \) can be viewed as a reduced simplicial set.
In this context,  the degeneracy operators and the face operators
 $d_0^{\PPP,n} $, $d_n^{\PPP,n} $
are  specified in (P3) and (P1), respectively. 
The other face operators are defined as
$d_i^{\PPP,n} = (\mathrm{id}_{{\PPP_{i-1}}} \times \nab_{2} \times \mathrm{id}_{{\PPP_{n-i-1}}})|_{\PPP_{n}}\,\, (0 < i < n)$.
This structure, in fact, forms a prepartial monoid with identity spine operators.
Furthermore, if $\PPP$ is a partial group, it becomes a prepartial group 
via $\nu_n (\pii) = \pii^{-1}$ for $\pii \in \PPP_n$.

\begin{lem}
Let $\PPP$ be a reduced simplicial set such that $\PPP_n \subseteq (\PPP_1)^n$ for $n \geq 1$,
and for each  $\pii = (\pi_1,\pi_2, \ldots, \pi_n) \in \PPP_n$,
we have the conditions 
 $d_0 (\pii ) =  (\pi_2, \ldots, \pi_n)$ and $d_n (\pii ) =  (\pi_1, \ldots, \pi_{n-1})$.
Then, $\PPP := ((\PPP_n)_n, (\nab_n)_n)$ forms a partial monoid with unit  element
$1 = s_0 (*)$, where 
$\nab_n = d_1^2 \circ \cdots \circ d_1^{\, n}: \PPP_n \to \PPP_1$
and $\PPP_0 = \{ * \}$.  
\end{lem}

\begin{proof}
For any simplicial set, the following holds:
\[
e_i^{\, p + 1 + r} \circ (d_{p+1}^{\, p + 2 + r} \circ \cdots \circ d_{p+1}^{\,n} ) = 
\begin{cases}
e_i^n & (1 \leq i \leq p)\\
(d_1^{\, 2} \circ \cdots \circ d_1^{\, q} )\circ (d_0^{\, p+q, p} \circ d_n^{\, n,r}) & (i=p+1)\\
e_{i+q-1}^n & (p+2 \leq i ),
\end{cases}
\]
where $n = p + q + r$ and $q > 0$.
By assumption,  for $\pii = (\pi_i)_i \in \PPP_n$, we have 
$(d_0^{\, p+q, p} \circ d_n^{\, n,r}) (\pii) = (\pi_{p+1},\ldots,\pi_{p+q})$,
and in particular, $e_i^n (\pii) = \pi_i$. 
Using these results, we derive the following:
\[
(\id_{\PPP_p} \times \nab_q \times \id_{\PPP_r}) (\pii)
=(d_{p+1}^{\,p+2+r} \circ \cdots \circ d_{p+1}^{\, n}) (\pii)
=:\rrho
\in \PPP_{p+1+r},
\]
%
\[
\nab_{p+1+r} \circ (\id_{\PPP_p} \times \nab_q \times \id_{\PPP_r}) (\pii)
= (d_{1}^{\,2} \circ \cdots \circ d_{1}^{\, p+1+r}) 
( \rrho ) 
= \nab_n (\pii).
\]
The proof of properties (P3) and (P4) follows straightforwardly.
\end{proof}

Now, let \( \SSS \) be a prepartial monoid. By the previous lemma, we conclude that 
$\eee \SSS := \left( \left( \eee^n(\SSS_n) \right)_n, \left( \nab_n^{\eee\SSS} \right)_n \right)$
is a partial monoid, 
where \(\nab_n^{\eee\SSS}\) is defined as the composition of the canonical isomorphism
$\eee^n(\SSS_n) \cong \SSS_n$
with the map
$d_1^2 \circ \cdots \circ d_1^n : \SSS_n \to \SSS_1 \; (= \eee_1(\SSS_1)).$
If \(\SSS\) is a prepartial group equipped with an inversion \(\nu\), then \(\eee\SSS\) becomes a partial group with inverse operation defined by 
$x \mapsto \nu(x)$ for $x \in \SSS_1$.
In this case, we call \(\eee\SSS\) the \textit{associated partial group} of \(\SSS\).

\section{Semidirect products}

In this section, we introduce the concept of the action of a partial group $\FFF$ 
on another partial group $\HHH$. 
We also define the corresponding semidirect product construction for partial groups.

Let $\HHH$ and $\FFF$ be  partial groups.
For each $r, n > 0$, 
let $\HHH_r \FFF_n$ be a designated subset of
$\HHH_r \times \FFF_n$ and let 
$\act_{r | n} =\act_{r | n}^{\HHH\FFF}: \HHH_r \FFF_n \to \HHH_r$,
$ (\eeta | \ggam) \mapsto \act_{r | n}  (\eeta | \ggam)  = \act (\eeta | \ggam)$
be a map.
We say that $\HHH := ( \HHH, ( \HHH_r \FFF_n )_{r,n}, (\act_{r | n})_{r,n} )$ 
is a \textit{(right)$\FFF$-partial group}  
if the following conditions hold:
\begin{description}[nosep]
\item[(AP1)]
    For each $r > 0$, 
    $(\HHH_r, (\HHH_r \FFF_n)_{n > 0}, (\act _{r | n})_{n>0} )$
    is an  $\FFF$-set.
\item[(AP2)]
For each \(r > 1\), the maps 
$(\eta_1, \ldots, \eta_r) \mapsto (\eta_2, \ldots, \eta_r)$
and
$(\eta_1, \ldots, \eta_r)$ $\mapsto (\eta_1, \ldots, \eta_{r-1})$
are maps of \(\FFF\)-sets from \(\HHH_r\) to \(\HHH_{r-1}\).
\item[(AP3)]
For each \( p, r \geq 0 \) and \( q > 0 \), the map 
$\bigl( \id_{\HHH_p} \times \nab_q \times \id_{\HHH_r} \bigr)\big|_{\HHH_{p+q+r}}$
is a map of \( \FFF \)-sets from
$\HHH_{p+q+r}$ to $\HHH_{p+1+r}$.
\item[(AP4)]
For each $0 \leq i \leq r$,
the map $(\eta_1, \ldots, \eta_r) \mapsto (\eta_1,\ldots \eta_i,1,$
$\eta_{i+1},\ldots, \eta_r)$
is a map of $\FFF$-sets from $\HHH_{ r}$ to $\HHH_{ r+1}$.  
\item[(AP5)]
    For each $\ggam \in \FFF_n$, we have that $(1 | \ggam ) \in \HHH_1 \FFF_n$
and $\act_{1 | n} (1 | \ggam ) = 1$.
\item[(AP6)]
   For each $r > 0$, the maps $\eeta \mapsto \eeta^{\pm 1} \times  \eeta^{\mp 1}$
   are maps of \( \FFF \)-sets from 
$\HHH_{ r }$ to $\HHH_{ 2r }$.
\end{description}

We note that the conditions (AP1)-(AP5) hold if and only if 
$\HHH$ is a simplicial object in the category of $\FFF$-sets.

\begin{ex}
\label{VactF}
(1)\,
Let $\HHH$ and $\FFF$ be partial groups,
and let $\psi : \FFF \to \mathrm{Aut} (\HHH)$ be a map of partial groups,
where $\mathrm{Aut} (\HHH)$ denotes the set of automorphisms of $\HHH$ 
with a (partial) group structure defined by 
$\eta (\phi_1 \phi_2) = (\eta \phi_1) \phi_2$ 
for
$\eta \in \HHH_1$ and $\phi_1, \phi_2 \in \mathrm{Aut} (\HHH)$.
Here, the action of \( \phi \in \operatorname{Aut}(\HHH) \) on \( \eta \in \HHH_1 \) 
is given by \( \eta\phi = \phi(\eta) \). 
Under this setup, \( \HHH \) becomes an \( \FFF \)-partial group 
with elements $\HHH_r \FFF_n := \HHH_r \times \FFF_n$
and operations defined by 
$\act _{r | n} (\eeta \,|\, \ggam):= (\eta_i \psi (\gam_1) \cdots \psi (\gam_n))_i$,
where $\eeta =(\eta_i)_i \in \HHH_r$ and $\ggam = (\gam_i)_i \in \FFF_n$.
In this case,  $\HHH$ is referred to as an $\FFF$-partial group of {\it Grazian-Henke type}
associated with $\psi$ 
 (\cite[Definition 3.1]{GH}). 
\\
(2)\,
Let $G$ be a group. Let $\FFF$ be a partial group and $X$ a friendly $\FFF$-set. 
Then the set ${\bigvee}_X G$ becomes an $\FFF$-partial group via the definition
\begin{multline*}
({\bigvee}_X G)_{r} \FFF_n
:=  \{
([ x, a_1 ] , [ x, a_2 ] , \ldots, [ x, a_r ] | \ggam ) \mid \\
( x| \ggam ) \in X \FFF_{ n },\, a_1, a_2, \ldots, a_r \in G \} ,
\end{multline*}
where the operation $\act$ is given by
\begin{multline*}
\act
( [ x, a_1 ] , [ x, a_2 ] , \ldots, [ x, a_r ] | \ggam)
:=\\
( [ \act (x| \ggam), a_1 ] , [ \act (x| \ggam), a_2 ] , \ldots, [ \act (x| \ggam), a_r ] ) .
\end{multline*}
Note that condition \((\mathrm{AP5})\) holds since \( X \) is \( \FFF \)-friendly. 
It is also notable that ${\bigvee}_X G$ is {\it not} an \( \FFF \)-partial group of 
Grazian-Henke type if $X\FFF_1 \subsetneq X \times \FFF_1$ and \( G \neq \{1\} \).
The former condition is satisfied, for example,  if $X = \{ 0, 1 \}$ and $\FFF = \Ups (\{ 0,1 \})$ 
(cf. Example \ref{UpP} (2) and Example \ref{PVP}).
\end{ex}

\begin{lem}
\label{LemFPG}
Let $\HHH$ be an $\FFF$-partial group.\\
 {\rm (1)}\, 
For each $ (\eta_1,\ldots, \eta_r | \ggam ) \in \HHH_r \FFF_n$,
we have $ (\eta_1 | \ggam ), \ldots,  (\eta_r | \ggam ) \in \HHH_1\FFF_n$, and
\begin{equation}
\act_{r | n} (\eta_1,\ldots, \eta_r | \ggam ) 
= 
(\act_{1 | n} (\eta_1 | \ggam ) , \ldots, \act_{1 | n} (\eta_r | \ggam )).
\end{equation}
 {\rm (2)}\, For each $ (\eta | \ggam ) \in \HHH_1 \FFF_n$,
we have $ (\eta^{-1} | \ggam ) \in \HHH_1\FFF_n$, and
\begin{equation}
\act_{1 | n} (\eta | \ggam )^{-1}
= 
\act_{1 | n} (\eta^{-1} | \ggam ).
\end{equation}
\end{lem}

\begin{proof}
(1)\,\, By (AP2) the projection 
$\HHH_r \to \HHH_1$, $(\eta_1, \ldots, \eta_r) \mapsto \eta_i$
is a map of \(\FFF\)-sets for each \(i\). 
Hence, applying this to \((\eta_1, \ldots, \eta_r | \ggam)\) yields the desired result.\\
(2) Let \((\eta | \ggam)$ be an element of $\HHH_1 \FFF_n\). By (AP6), we have 
$(\eta, \eta^{-1} | \ggam) \in \HHH_2 \FFF_n$
and
$\act(\eta, \eta^{-1} | \ggam) = ( \act(\eta |\ggam),\, \act(\eta | \ggam)^{-1} ).$
On the other hand, by Part (1),
$\act(\eta, \eta^{-1} | \ggam)$ $=$ $( \act(\eta | \ggam),\, \act(\eta^{-1} | \ggam) ).$
Comparing the second components, we obtain the desired equality in (2). 
\end{proof}

Let $\HHH$ be an $\FFF$-partial group.
It is straightforward to verify that the sequence $\FHs:= (\HHH_n \FFF_n)_n$  of sets
forms a reduced simplicial set with the following operations:
\begin{equation}
d_i  (\xxi | \ddel ) 
=
\begin{cases}
(d_i (\xxi) | d_{n-i}(\ddel) ) & (i < n)\\
(d_n (\act_{n | 1} (\xxi | \del_1)) | d_{0}(\ddel) ) & (i = n), \\
\end{cases}
\end{equation}
\begin{equation}
s_i  (\xxi | \ddel ) 
=
(s_i (\xxi) | s_{n-i}(\ddel) ), 
\end{equation}
for every  $(\xxi | \ddel) = (\xi_1,\ldots, \xi_n|\del_1,\ldots,\del_n) \in \HHH_n \FFF_n$.
We refer to  $\FHs$ as the {\it separated semidirect product} of $\HHH$ by $\FFF$.

To show that $\FHs$ is a prepartial group, 
it is convenient to identify an element \((\xi | \delta) \in \HHH_1 \FFF_1\) with
$\delta^{-1} \ltimes \xi := (\delta^{-1}, \xi) \in \FFF_1 \times \HHH_1$.
We then introduce the following sets:
\[
(\FH)_1 = (\FH)_1^\prime
 := \{ \delta^{-1} \ltimes \xi \mid (\xi | \delta) \in \HHH_1 \FFF_1 \},
\]
and for \( n > 1 \),
\[
(\FH)_n^\prime
:=
\{
 (\gam_i \lt \eta_i)_i \in ( \FFF \lt \HHH)_1^{\, n} \mid 
 (\eta_j | \gam_{j+1}, \ldots, \gam_n) \in \HHH_1 \FFF_{n-j}\,\,
 (1 \leq j < n)
\}.
\]
Next, define a map \( A \) by
\[
 A ( \aalp )
 :=
( \act (\eta_1 | \gam_2, \ldots, \gam_n) , \act (\eta_2 | \gam_3, \ldots, \gam_n) ,
\ldots,
\act (\eta_{n-1} | \gam_n) , \eta_n ),
\]
and extend it to a map \( A^+ \) as
\[
 A^+ (\aalp) := ( A ( \aalp ) | \gam_n^{-1}, \ldots, \gam_1^{-1} ),
\]
where $\aalp =  (\gam_i \lt \eta_i)_i \in (\FH)_n^\prime$. 
\begin{lem}
For each $(\xxi | \ddel ) = ((\xi_i)_i | (\del_i )_i) \in ( \FHs )_n$, 
we have $\aalp = (\gam_i \lt \eta_i)_i := \eee^n (\xxi | \ddel ) \in ( \FH )_n^\prime$.
Furthermore,  $A^+ (\aalp ) = (\xxi |\ddel )$. 
\end{lem}

\begin{proof}
For each $1 \leq i \leq n$, we have:
\begin{equation}
\label{exidel}
e^n_i (\xxi | \ddel )
=
(\act (\xi_i | \del_1,\ldots, \del_{n-i}) | \del_{n+1-i} ).
\end{equation}
By conditions (AP2), (A1), (A6), and (A2), this implies:
\begin{equation}
\label{eta|gam}
(\eta_i | \gam_{i+1},\ldots, \gam_n ) = (\act (\xi_i | \del_1, \ldots, \del_{n-i})| \del_{n-i}^{-1},\ldots, \del_{1}^{-1}) \in \HHH_1\FFF_{n-i}.
\end{equation}
The second assertion follows directly from \eqref{eta|gam}.
\end{proof}

\begin{prop}
For each $\FFF$-partial group $\HHH$,  $\FHs$ is a prepartial group with 
inversion $\nu_n (\xxi | \ddel ) = (\act_{n|n} (\xxi^{-1} | \ddel ) | \ddel^{-1} )$,  
where $(\xxi | \ddel ) \in \FFF_n \HHH_n$.
The associated partial group $\FH := \eee ( \FHs )$ 
is given as follows:
\begin{align}
\label{FHn}
( \FH )_n
&:=
\{
\aalp \in ( \FH )_n^\prime \mid 
A^+ ( \aalp ) \in \HHH_n \FFF_n
\},\\
\label{nabFH}
\nab^{\FFF \ltimes \HHH}_n ( \aalp )
& := \nab_n^{\FFF} ((\gam_i)_i) \lt \nab_n^{\HHH} ( A ( \aalp ))
\quad (\aalp = (\gam_i \lt \eta_i)_i \in ( \FH )_n),\\
\label{1FH}
1_{\FFF \ltimes \HHH} 
& = 1_\FFF \lt  1_\HHH,\\
\label{gameta-}
( \gam \lt \eta )^{-1}
& =
 \gam^{-1} \lt \act ( \eta | \gam^{-1} )^{-1}
 = \gam^{-1} \lt \act ( \eta^{-1} | \gam^{-1} ).
\end{align}
\end{prop}

\begin{proof}
By the above lemma, \(\FHs\) is a prepartial monoid and 
\(\eee(\FHs)_n\) is a subset of \((\FH)_n\) as defined in \eqref{FHn}.  
Let $\ssig = (\xxi | \ddel )$ be an element of $\HHH_n \FFF_n$.
By conditions (AP6) and (A6), we have
$\ttau:= (\xxi^{-1} \times \xxi | \ddel \times \ddel^{-1} )\in \HHH_{2n} \FFF_{2n}$.
Therefore, by condition (AP2), $\nu_n (\ssig)$ is a well-defined element of $\HHH_n\FFF_n$.
It is easy to verify that $\nu_n$ satisfies condition (Inv1) and that $\ttau$ satisfies condition (Inv2).
Thus, $\FHs$ is a prepartial group.
Let
$\aalp = (\gam_i \lt \eta_i)_i$ be an element of $( \FH )_n$.
Since 
$(\eta_i | \gam_{i+1},\ldots,\gam_n,$ $\gam_n^{-1},\ldots,\gam_{i+1}^{-1}) \in \HHH_1 \FFF_{2(n-i)}$,
we have
$\eta_i = \act ( \xi_i | \del_1, \ldots, \del_{n-i})$,
where $((\xi_i)_i | (\del_i )_i) = A^+ ( \aalp )$. 
Hence, by \eqref{exidel}, we obtain
$\eee ( \FHs )_n \ni \eee^n ((\xi_i)_i | (\del_i )_i) = \aalp$.
This shows that $\eee ( \FHs )_n = ( \FH )_n$.
The formula \eqref{nabFH} follows from 
\[
(d_1^{\, 2} \circ \cdots \circ  d_1^{\, n} ) ( \xxi | \ddel )
=
\nab_n (\ddel )^{-1} \lt \nab_n ( \xxi).
\]
The proofs of \eqref{1FH} and \eqref{gameta-} are straightforward.
\end{proof}

We refer to $\FFF \lt \HHH$ as the {\it semidirect product} of $\HHH$ by $\FFF$.

\begin{ex}
\label{GHSD}
Let $\HHH$ be the  $\FFF$-partial group of Grazian-Henke type
associated with $\psi\,: \FFF \to \mathrm{Aut} (\HHH)$. 
By \eqref{FHn}, we have
\begin{equation}
\label{GHSDn} 
 (\FFF \ltimes \HHH)_n
 =
\bigl\{ 
\aalp = (( \gam_i, \eta_i ))_i \in (\FFF_1 \times \HHH_1)^n
\mid 
( \gam_i )_i \in \FFF_n, \,\, 
A ( \aalp ) \in 
\HHH_n 
\bigr\}.
\end{equation}
This demonstrates that our semidirect product construction extends 
the framework presented in  \cite{GH}. 
\end{ex}

\section{A characterization of the semidirect product}

In Section 10, a semidirect product $\FH$ is defined in terms of simplicial theory.
In this section, we characterize $\FH$ by means of its universal mapping property.
This result will serve as a foundation for proving the theorems presented in the following section.

Let $\HHH$ and $\PPP$ be $\FFF$-partial groups. A map $f: \HHH \to \PPP$ 
of partial groups is called a {\it map of $\FFF$-partial groups} if  
$f^n |_{\HHH_n}: \HHH_n \to \PPP_n$ is a map of $\FFF$-sets, for each $n>0$.
We denote by $\bf{PG}_\FFF$ the resulting category 
of $\FFF$-partial groups and maps of $\FFF$-partial groups.
To begin, we construct a distinguished object $\PPP^{\rm ad}$ in $\bf{PG}_\PPP$,
referred to as the {\it adjoint} (or {\it conjugate}) $\PPP$-{\it partial group} 
associated with an  arbitrary partial group $\PPP$ .

\begin{lem}
\label{PadP}
 Let $\PPP$ be a partial group.
 Define $\PPP^{\rm ad}_{r} \PPP_n \subseteq \PPP_r \times \PPP_n$ and
 $ \act^{{\rm ad}}_{\,r | n}: \PPP^{\rm ad}_{r} \PPP_n \to \PPP_r$ by
\begin{align*}
\PPP^{\rm ad}_{r} \PPP_n
 &:=
  \{ ( \pii | \ggam ) \in \PPP_r \times \PPP_n \mid 
\rrho ( \pii | \ggam )
  \in \PPP_{(2n+1)r} \},\\
   \act^{{\rm ad}}_{\,r | n}  ( \pii | \ggam )
& : = 
  ( \nab_{2n+1} (\rrho_i (\pii | \ggam )))_{1 \leq i \leq r} ,
\end{align*}
where
\begin{align*}
  \rrho ( \pii | \ggam )
  &: = \rrho_1 ( \pii  | \ggam ) \times \cdots \times \rrho_r ( \pii | \ggam ),\\
  \rrho_j ( (\pi_i)_i | \ggam )
  &:=  
  \ggam^{-1} \times ( \pi_j ) \times \ggam
  \quad ( 1 \leq j \leq r).
\end{align*}
Then $\PPP^{\rm ad} := ( \PPP, ( \PPP^{\rm ad}_{r} \PPP_n )_{r,n}, (\act^{{\rm ad}}_{\,r | n})_{r,n} )$ 
is a $\PPP$-partial group.
\end{lem}

\begin{proof}
The proof is elementary. 
We only provide the proof for (A2) and (AP3) here.
Let $( \pii | \ggam )$ be an element of $\PPP^{\rm ad}_{r} \PPP_n$.
For the proof of 
(A2),  set
$( \pii^\pr | \ggam^\pr ) := ( \act^\ad_{r | p} \times \id_{\PPP_q} )( \pii | \ggam )$,
where $n = p + q$.
Then we have 
\begin{equation*}
\rrho_j ( \pii^\pr | \ggam^\pr ) 
=
( {\id_{{\PPP_q} }} \times \nab_{2p +1} \times {\id_{{\PPP_q} }}  )
( \rrho_j (\pii | \ggam )).
\end{equation*}
Hence  $( \pii^\pr | \ggam^\pr ) \in \PPP^{\rm ad}_{r} \PPP_q$ by \eqref{P2-1}.
 We also obtain
\begin{align*}
\act^\ad_{\, r | q} 
( \pii^\pr | \ggam^\pr ) 
& =
\left(
\nab_{2q+1} \circ ( {\id_{{\PPP_q}}}  \times \nab_{2p +1} \times {\id_{{\PPP_q} }} )
(\rrho_j ( \pii | \ggam ))
\right)_j\\
& =
\act^\ad_{\, r | p+q} (\pii | \ggam ),
\end{align*}
showing that (A2) holds.
For (AP3), set
$\pii^\prr := ( {\id_{{\PPP_s}}} \times \nab_{t} \times {\id_{{\PPP_u} }} ) ( \pii )$,
where  $r = s + t +u$.
Then,
\begin{gather*}
\rrho ( \pii^\prr | \ggam )
 =
( \id_{{\PPP_{\! 2ns+n+s} }} \times \nab_{ 2nt -2n+t } \times \id_{{\PPP_{\! 2nu+n+u} }}  )
( \rrho (\pii | \ggam )) \in \PPP_{(2n+1)(s+1+u)}.
\end{gather*}
Moreover,
\begin{gather*}
({\id_{{\PPP_s}}}  \times \nab_t \times {\id_{{\PPP_u}}}  ) ( \act^\ad_{\, r | n} (\pii | \ggam ))\\
=
(\underbrace{(\nab_{2n+1}\times \cdots \times  \nab_{2n+1} )}_s
\times \nab_{(2n+1)t} \times 
\underbrace{(\nab_{2n+1}\times \cdots \times  \nab_{2n+1} )}_u)
(\rrho ( \pii | \ggam ))\\
 =
\act^\ad_{s+1+u | n} (\pii^\prr | \ggam ).
\end{gather*}
Thus, we have proved (AP3).
\end{proof}

Let $\phi:\FFF \to \HHH$ be a map of partial groups,
and let $\PPP$ be an $\HHH$-partial group.
It is easy to see that $\PPP$ becomes an $\FFF$-partial group via
\begin{gather}
\label{KF}
\PPP_r \FFF_n
:= \{ ( \pii | \ggam ) \in  \PPP_r \times \FFF_n
\mid  ( \pii | \phi^n (\ggam )) \in \PPP^{}_r \HHH_n \},\\
\label{actKF}
\act_{r | n}^{\PPP\FFF} (\pii | \ggam ):=  \act_{r | n}^{\PPP\HHH}  (\pii | \phi^n (\ggam ))
\quad ( (\pii | \ggam ) \in  \PPP_r \FFF_n ).
\end{gather}
In particular, $\HHH$ itself becomes an $\FFF$-partial group $\HHH^\phi$ via   
\begin{gather}
\HHH^\phi_r \FFF_n
:= \{ ( \eeta | \ggam ) \in  \HHH_r \times \FFF_n \mid 
( \eeta | \phi^n (\ggam )) \in \HHH^{\rm ad}_r \HHH_n \},\\
\act_{r | n}^\phi ( \eeta | \ggam ):=  \act^{{\rm ad}}_{r | n}  ( \eeta | \phi^n (\ggam ))
\quad (  ( \eeta | \ggam ) \in  \HHH^\phi_r \FFF_n ). 
\end{gather}

The assignment $\phi \mapsto \HHH^\phi$ defines a  
functor $\Phi$ from the comma category
$( \FFF\downarrow\bf{PG} )$ to $\bf{PG}_\FFF$. 
Here, the objects of $( \FFF\downarrow\bf{PG} )$
are maps $\FFF \stackrel{\phi}{\to} \HHH$ of partial groups,
 and morphisms
from  $\FFF \stackrel{\phi}{\to}  \HHH$ to $\FFF \stackrel{\psi}{\to}  \PPP$ are maps $\chi:\HHH \to \PPP$
of partial groups such that $\chi \circ \phi = \psi$.
Now we will show that our semidirect product gives a left adjoint of the functor $\Phi$.
\begin{thm}
\label{UMP}
For each partial group $\FFF$ and each $\FFF$-partial group $\HHH$, 
there exists a universal arrow $(\FFF \stackrel{\iota}{\to} \FH, h)$
from $\HHH$ to the functor $\Phi$.
In other words, there exist
a map $\FFF \stackrel{\iota}{\to} \FH$ of partial groups
and a map $h:\HHH \to ( \FH )^\iota$ of $\FFF$-partial groups satisfying the following property:
for each map $\FFF  \stackrel{\phi}{\to}  \PPP$ of partial groups 
and each map $f: \HHH \to \PPP^\phi$ of $\FFF$-partial groups,
there exists a unique map $\psi :  \FH \to \PPP$ of partial groups such that
the following diagrams commute:

\begin{equation}
\label{CD}
   \xymatrix{
    & \FFF \ar[ld]_\iota  \ar[rd]^\phi & \\
    \FH     \ar[rr]_{\psi} && \PPP
   },
   \quad
   \xymatrix{
\HHH \ar[dr]_f  \ar[rr]^{h}  & & \FH   \ar[dl]^{\psi} \\
& \PPP  &
}.
\end{equation}
\bigskip\par\noindent
The maps  $\iota$, $h$, and $\psi$ are explicitly given by the following formulas:
\begin{equation}
\label{ihpsi}
\iota ( \gam ) =  \gam \lt 1,\quad
h( \eta ) = 1 \lt \eta,\quad
\psi ( \gam \lt \eta ) = \nab_2 ( \phi ( \gam  ), f( \eta )),
\end{equation}
where $\gam \in \FFF_1$ and $\eta \in \HHH_1$.
\end{thm}

\begin{proof}
Let  $\iota$, $h$, and $\psi$ be as defined in \eqref{ihpsi}.
It is straightforward to verify that the map $\iota$ is well-defined as a map of partial groups, 
and that $h$ is a map of $\FFF$-partial groups.
We now show that $\psi$ is also a well-defined map of partial groups.

Let $\aalp = ( \gam_i \lt \eta_i )_i$  be an element of $( \FH )_n$.
Since $A^+ ( \aalp ) \in \HHH_n \FFF_n$ and $f$ is a map of $\FFF$-sets, we deduce that
\[
( \act_{1 | n}^\phi ( f( \eta_1 ) | \gam_2, \ldots,  \gam_n ), \ldots, 
\act_{1 | 1}^\phi ( f( \eta_{n-1} ) | \gam_n ),f( \eta_n ) | \ggam^{-1} )
\in
\PPP^\phi_n \FFF_n,
\]
where $\ggam =( \gam_i )_i$.
Equivalently, we have $\bbeta \in \PPP_{ n(2n+1)}$, where
\begin{gather*}
\bbeta := \phi^n ( \ggam ) \times 
( \act^{{\rm ad}}_{1 | n}  ( f( \eta_1 ) | \phi ( \gam_2 ), \ldots,  \phi ( \gam_n )))
\times \phi^n ( \ggam )^{-1}  \times  \cdots \\
\cdots \times  \phi^n ( \ggam ) \times 
(\act^{{\rm ad}}_{1 | 1}  ( f( \eta_{n-1} ) | \phi ( \gam_n ))) 
\times \phi^n ( \ggam )^{-1} \times  \phi^n ( \ggam ) \times (f( \eta_n )) \times \phi^n ( \ggam)^{-1} .
\end{gather*}
From this, it is easy to deduce that
\[
( \phi ( \gam_1  ), f( \eta_1 ), \ldots, \phi ( \gam_n  ), f( \eta_n )) \in \PPP_{2 n}.
\]
Hence $\psi$ is well-defined and satisfies 
$\psi^n ( \aalp ) \in \PPP_n$.
Moreover, we have
\[
\psi (\nab ( \aalp ))
=
\nab (\bbeta \times \phi^n ( \ggam ))
=
\nab (\psi^n ( \aalp )).
\] 
Thus, we obtain a map $\psi : \FH \to \PPP$ of partial groups 
which makes the diagrams in the theorem commute.
Finally, suppose that $\psi' : \FFF \lt \HHH \to \PPP$ is another map 
satisfying the same conditions as $\psi$, and let $\gam \lt \eta \in (\FH)_1$. 
Since  $(\gam \lt 1, 1\lt \eta) \in ( \FH )_2$ and 
$\nab_2(\gam \lt  1, 1 \lt \eta) = \gam \lt \eta$,
we have
$\psi^\prime (\gam \lt \eta ) = \nab_2 (\psi^\prime (\gam \lt  1), \psi^\prime (1 \lt  \eta )) = \psi (\gam \lt  \eta)$.
This completes the proof of Theorem \ref{UMP}. 
\end{proof}

\section{Figures and semidirect products}
In this section, we relate our semidirect product construction to certain classes of figures. We begin by considering figures with two congruent connected components that do not satisfy the condition in Theorem \ref{P=C2GG}, 
as illustrated in the
right panel of Figure 5.

\begin{thm}
\label{P(F0F1)}
Let \( F \subseteq \mathbb{R}^d \) be a figure 
consisting of two congruent connected components, denoted by $F_{+}$ and $F_{-}$. 
Suppose there is no isometry \( a \in E(d) \) satisfying \( F_{\pm} a = F_{\mp} \) simultaneously. 
Then we have
\[
 \PPP  (F)  \cong 
\Ups (\{0,1\}) \lt \left( G(F_+) \vee_{G(F_+) \cap G(F_-)} G(F_-) \right),
\]
%
where $\Ups ( \{0, 1 \} )$ is defined as in \eqref{V01n}. If $F_+ a = F_-$ for some $a \in M(d)$,
then we also have 
\begin{align*}
 S \PPP  (F) & \cong 
\Ups (\{0,1\}) \lt \left( SG(F_+) \vee_{SG(F_+) \cap SG(F_-)} SG(F_-) \right).
\end{align*}
\end{thm}

The proof of the above theorem is very similar to that of the subsequent theorem; therefore, we present only the latter.

\begin{thm}
\label{SimP=}
Let \( F \subseteq \mathbb{R}^d \) consist of two similar connected components \( F_{\pm} \). 
Suppose that $F_+$ and $F_-$ are not congruent
and $0 < d( F_\pm ) < \infty$.
Then, we have
\[
 Sim\PPP  (F) \cong 
\Ups (\{0,1\}) \lt \left( G(F_+) \vee_{G(F_+) \cap G(F_-)} G(F_-) \right).
\]
\end{thm}

\begin{proof}
We use the following notations:
$F_0 := F_-$; $F_1 :=F_+$; $G_0 := G(F_0)$; $G_1 := G(F_1)$; 
$\PPP:= Sim \PPP (F)$; $\FFF:= \Ups ( \{0, 1 \})$; 
$\HHH := G_0 \cup G_1 \cong G_0 \vee_{G_0 \cap G_1} G_1 $. 
Let $\gam$ be an element of $Sim(d)$ such that  $F_0 \gam = F_1$.
Then we have
\begin{gather}
\label{Fg=F}
F_{\iota_0} \gam^{\iota_1 - \iota_0} = F_{\iota_1},
\quad
\gam^{\iota_0 - \iota_1} G_{\iota_0}
\gam^{\iota_1 - \iota_0} =  G_{\iota_1}
\quad
(\iota_0, \iota_1 \in \{0, 1 \} ),\\
\label{K1=}
\PPP_1 = \HHH_1  \amalg \gam G_1  \amalg  \gam^{-1} G_0.
\end{gather}
%
%
Here we have derived  \eqref{K1=} from\eqref{Pxy}  and Proposition \ref{SimG=G}.
Using \eqref{V01n} and \eqref{Fg=F},
we see that \( \delta \mapsto \gamma^\delta \) (where \( \delta = 0, \pm 1 \)) gives 
a map $\phi :  \FFF \to \PPP$ of partial groups.
Also, using \eqref{Fg=F}, we see that  $\HHH$ becomes an $\FFF$-partial group via
\begin{gather}
\label{HrFn}
\HHH_r \FFF_n 
= 
\left\{ ( \aaa | \iota_1 - \iota_0 \ldots,  \iota_n - \iota_{n-1} ) \mid 
\iota_0, \ldots, \iota_n \in \{ 0, 1 \}, 
\aaa \in  {G_{\iota_0}}^{\! r}  \right\},\\
\act_{r | n}^{\HHH\FFF}  ( ( a_i )_i | \iota_1 - \iota_0 \ldots,  \iota_n - \iota_{n-1} )
= (\gam^{\iota_0 - \iota_n} a_i \gam^{\iota_n - \iota_0} )_i.
\end{gather}
To show that the inclusion $f\!:\HHH \hookrightarrow \PPP$ is a map of $\FFF$-partial groups,
choose
$\iota_0,\ldots, \iota_n \in \{ 0, 1 \}$ and
$a_1,\ldots, a_r \in G_{\iota_0}$.
By \eqref{Fg=F}, 
$\ggam^{-1} \times ( a_1 ) \times \ggam \times \cdots \times \ggam^{-1} \times ( a_r ) \times \ggam$ 
is a parade with a starting point $F_{\iota_n}$, where
$\ggam := (\gam^{\iota_1 - \iota_0}, \ldots, \gam^{\iota_n - \iota_{n-1}})$.
Hence, we have $( ( a_i )_i | \ggam) \in \PPP^\ad_r \PPP_n$, or
$( ( a_i )_i | (\iota_i - \iota_{i-1} )_i) \in \PPP^\phi_r \FFF_n$.
Moreover, we have
\begin{gather*}
\act_{r | n}^\phi (( a_i )_i | (\iota_i - \iota_{i-1} )_i) 
= 
\act_{r | n}^\ad (( a_i )_i | \ggam) \\
=
(\gam^{\iota_0 - \iota_{n}} a_i \gam^{\iota_n - \iota_{0}})_i
=
\act_{r | n}^{\HHH\FFF} (( a_i )_i | (\iota_i - \iota_{i-1} )_i)
\end{gather*}
as desired.
Applying Theorem \ref{UMP}, we obtain a map $\psi: \FH \to \PPP$ of partial groups,
which sends  $\del \lt \eta \in ( \FH )_1$ to $ \nab_2(\gam^\del , \eta)$.
By \eqref{K1=} and \eqref{HrFn}, $\psi$ gives a bijection $( \FH )_1 \to \PPP_1$. 
It remains to show that $\PPP_n \subseteq \psi^n ( ( \FH )_n )$.
Let $\pii = (\pi_i )_i$ be an element of $\PPP_n$,
which defines a route $F_{\iota_0},\ldots, F_{\iota_n}$.  
Since $\eta_i := \gam^{\iota_{i-1} - \iota_{i}} \pi_i \in G_{\iota_i }$,
we see that $\xi_i := ( \iota_i - \iota_{i-1} ) \lt  \eta_i $  is an element of $( \FH )_1$.
Since $\psi (\xi_i ) = \pi_i$, 
it suffices to show that
 $( \FH )_n$ contains: 
 \begin{align}
(\FH )_n^{\prime \prime }
& :=
\left\{ 
( (\iota_1 - \iota_{0} ) \lt \eta_1, \ldots, (\iota_n - \iota_{n-1} ) \lt \eta_n ) \mid 
\right. 
\nonumber
\\
& \left.
\iota_0, \ldots, \iota_n \in \{ 0, 1 \}, 
\eta_1 \in G_{\iota_1}, \ldots, \eta_n \in G_{\iota_n} 
\right\}.
\end{align}
Let  $\aalp =  ( (\iota_i - \iota_{i-1} ) \lt  \eta_i ))_i$ be an element
of $( \FH )^{\prime \prime }_n$.
By \eqref{Fg=F}, we have:
\[
A^+ (\aalp) 
= 
(\gam^{\iota_1 - \iota_n} \eta_1 \gam^{\iota_n - \iota_1},
\gam^{\iota_2 - \iota_n} \eta_2 \gam^{\iota_n - \iota_2}, \ldots, \eta_n
| \iota_{n-1} - \iota_n, \ldots, \iota_{0} - \iota_{1})
\in \HHH_n \FFF_{n}.
\]
This completes the proof of Theorem \ref{SimP=}.
\end{proof}



Finally, we will consider a class of figures that includes the Six Coins as an example.

\begin{thm}
\label{P(F+x)}
Let $F$ be a subset of ${\Bbb R}^d$ such that 
$cc (F) = \{ F_0 + x\,|\, x \in P \}$
for some $F_0, P \subset \R^d$. 
Assume that $G(F_0)$ is a subgroup of $O(d)$.
Then, we have the following isomorphisms:
\begin{gather}
\label{PF=}
 \PPP (F)  \cong 
\left( \Ups ({P}) \lt {\bigvee}_{P} G( F_0 ) \right) /\!\!\sim,\\
\label{SPF=}
S \PPP (F)  \cong 
\left( \Ups ({P}) \lt {\bigvee}_{P} SG( F_0 ) \right) /\!\!\sim
\end{gather}
for some congruence $\sim$,
where the $\Ups (P)$-partial group structures of ${\bigvee}_{P} G( F_0 )$
and $ {\bigvee}_{P} SG( F_0 )$ are as in Example \ref{VactF} (2).
Specifically, the congruence \( \sim \) is defined by:
 \begin{equation}
\label{sim}
(y - x ) \lt  [ y, a ] \sim (y^\pr - x^\pr ) \lt [ y^\pr, a^\pr ]
\Leftrightarrow
a = a^\pr,\,
(x^\pr - x) a = y^\pr - y.
\end{equation}
\end{thm}
%

\begin{proof}
We construct a map $\psi : \FH \to \PPP$ of
partial groups using Theorem \ref{UMP},
where $\FFF:= \Ups (P)$, $\HHH :=  {\bigvee}_{P} G( F_0 )$, and 
$ \PPP := \PPP (F)$.
Clearly, the map $y - x \mapsto t_{y-x}$\, (where $x, y \in P$) gives a map
$\phi : \FFF \to \PPP$  of partial groups.
Applying  the universal mapping property of the wedge sum to
the inclusions $G( F_0 + x ) \hookrightarrow \PPP$ (for \( x \in P \)),
we obtain a map  $f : \HHH \to \PPP$ 
of partial groups, given by $f ( [x, a] ) = t_ x^{\, -1} a t_x$. 
To verify that $f$ is a map $\HHH \to \PPP^\phi$ of $\FFF$-partial groups,
let 
$( \eeta | \ddel ) = (([x_0,a_i ])_i | x_1-x_0,\ldots, x_n - x_{n-1} )$
be an element of $\HHH_r \FFF_n$,
and let $\ggam \in \PPP_n$ be $\phi^n (\ddel)$.
Since
$\ggam^{-1} \times (t_{x_0}^{\, -1} a_1 t_{x_0}) \times \ggam \times \cdots \times
\ggam^{-1} \times (t_{x_0}^{\, -1} a_r t_{x_0}) \times \ggam$ 
is a parade with a starting point $F_0 + x_n$,
we conclude that $( f^r ( \eeta ) | \ddel )$
is an element of $\PPP^\phi_r \FFF_n$.
Moreover, we have:
\[
\act^\phi_{r | n} ( f^r ( \eeta ) | \ddel ) 
=
( \nab ( \ggam^{-1} \times (t_{x_0 }^{\, -1} a_i t_{x_0}) \times \ggam ))_i
=
( t_{x_n }^{\, -1} a_i t_{x_n} )_i
= f^r (\act_{r | n}^{\HHH\FFF} (\eeta | \ddel ))
\]
as desired.
By applying Proposition \ref{FundT} (2) to the obtained map
$\psi$, we establish an isomorphism 
$ ( \FH )/ \!\!\sim_\psi \, \cong {\rm Im} (\psi )$
of partial groups.
%
Let $\pii = (\pi_i)_i$ be an element of $\PPP_n$, 
and suppose it defines a route $F_0 + x_0, F_0 + x_1,\ldots,F_0 + x_n$.
Then, we have
$a_i := t_{x_{i-1}} \pi_i t_{x_i}^{\, -1} \in G(F_0 )$ for each $1 \leq i \leq n$.
It is straightforward to verify that 
$\aalp: = (( x_i - x_{i-1} ) \lt [x_i, a_i] )_i$
is an element of 
$ ( \FH )_n$ and that
$\psi^n ( \aalp ) = \pii$.
This shows that ${\rm Im } ( \psi )= \PPP$, completing the proof of \eqref{PF=}.
The proof of
\eqref{SPF=} follows a similar approach. 
Finally, suppose that $\psi ((y - x ) \lt  [ y, a ] ) = \psi ((y^\pr - x^\pr ) \lt [ y^\pr, a^\pr ])$ for
$x, y,  x^\pr, y^\pr \in P$ and $a, a^\pr \in G(F_0)$. 
Then, we have
$a t_{(x^\pr - x) a} = t_{x^\pr - x} a = a^\pr t_{y^\pr - y}$.
Since $G(F_0) \leq O(d)$, this is equivalent to the right-hand side of \eqref{sim}.
\end{proof}

Let $F$ be as in the above theorem. We say that $F$ is of 
{\it SDP (semidirect product) type} if the theorem provides an isomorphism
$ \PPP (F) \cong \Ups ({P}) \lt {\bigvee}_{P} G( F_0 )$.
Roughly speaking, $F$ is of SDP type 
when the only non-translation symmetries in $F$ are those of $F_0$.

\begin{prop}
\label{sdcond}
{\rm  (1)}\,
If $|P| \geq 2$ and $G(F) \not\subset \{ t_x\,|\, x \in \Bbb{R}^d \}$, then $F$ is not of SDP type.\\
{\rm  (2)}\,
Suppose $|v| = |w|$ for $v, w \in \Ups ({P})_1$ implies $v = \pm w$.
Also, suppose $va = \pm v$ for $v \in \Ups ({P})_1$ and $a \in G(F_0)$ only if 
either $v = 0$ or $a = 1$. Then $F$ is of SDP type.
\end{prop}

\begin{proof}
(1)\,
Let $y_1$ and $y_2$ be two distinct elements of $P$ and 
let  $g \in G(F)$ be an element that is not a translation.
For $i = 1,2$, define $x_i \in P$, $a_i \in G(F_0)$, and 
$\alp_i \in (\Ups ({P}) \lt {\bigvee}_{P} G( F_0 ))_1$ by 
$F_0 + x_i = ( F_0 + y_i ) g^{-1}$, $a_i = t_{x_i} g t_{- y_i}$, and
$\alp_i = (y_i - x_i)  \lt [ y_i, a_i ]$, respectively. 
Since $a_1, a_2 \ne 1$,  we have $\alp_1 \ne \alp_2$.
Hence, Part  (1) follows from $\psi ( \alp_1 ) = g = \psi ( \alp_2 ) $. \\
(2)\,
Let $\alp_1 = (y_1- x_1)  \lt [ y_1, a_1 ]$ 
and $\alp_2 = (y_2- x_2)  \lt [ y_2, a_2 ]$ be elements of 
$(\Ups ({P}) \lt {\bigvee}_{P} G( F_0 ))_1$ such that 
$\psi ( \alp_1 ) = \psi ( \alp_2 ) $.
By \eqref{sim} and the first assumption, we have $ x_2 - x_1 = \epsilon ( y_2 - y_1 )$,
where $\epsilon = \pm 1$.
When $x_1 = x_2$, this obviously implies $\alp_1 = \alp_2$.
When $x_1 \ne x_2$, the second assumption implies $a_1 = 1$.
Hence, we have $\alp_1 = \alp_2$. 
\end{proof}

\begin{figure}[h]
\centering
    \includegraphics[width=0.25\linewidth]{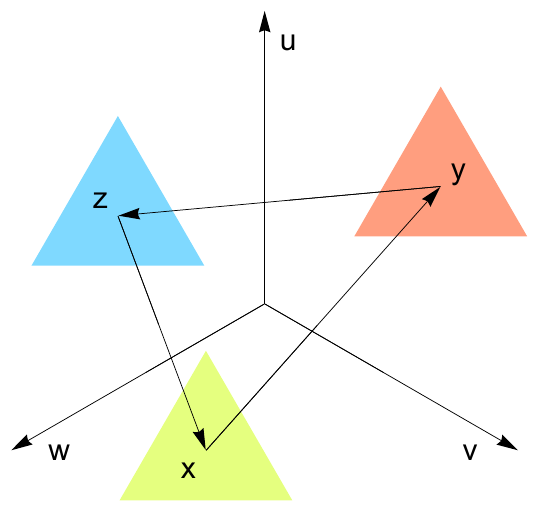}
\caption{A figure of SDP type.}\label{fig7}
\end{figure}

\begin{ex}
(1)\,
By Part (1) of the above Proposition, the figure of the Six Coins is not of SDP type. \\
(2)\,
Let $P = \{ x, y, z \}  \subset \R^2$ and $F$ be as in Figure 7. Since 
$| y - x | > | z - y | >| x - z |$, and
none of the vectors $y - x$, $z - y$, and $x - z$ are perpendicular to 
$u=(0,1)$, $v = ( \sqrt{3}/2,-1/2 )$, or $w =  ( -\sqrt{3}/2,-1/2 )$,
it follows from Part (2) of the above Proposition that $F$ is of SDP type. 
That is, we have 
$ \PPP (F) \cong \Ups ( \{ x, y, z \}) \lt  \left( D_3 \vee D_3  \vee D_3 \right)$.
\end{ex}


\bigskip\bigskip
\noindent
Graduate School of Mathematics, Nagoya University, Chikusa-ku, Nagoya, 464-8604, Japan\\
{\it{Email address}}: \texttt{hayashi@math.nagoya-u.ac.jp}


\end{document}